\documentclass[12pt,leqno]{article}
\usepackage{amssymb,url}

\newtheorem{proposition}{Proposition}[section]
\newtheorem{lemma}[proposition]{Lemma}

\newtheorem{theorem}[proposition]{Theorem}

\newtheorem{tmpdefinition}[proposition]{Definition}

\newtheorem{tmpobservation}[proposition]{Remark}
\newenvironment{observation}{\begin{tmpobservation}\rm }{\end{tmpobservation}}

\newenvironment{proof}{\par\topsep6pt plus 6pt\relax
\rm\trivlist\item[\hskip\labelsep \it Proof. ]}{\endtrivlist}
\newcommand{\cqd}{\hfill $\Box$\par}
\newcommand{\qedhere}{\mbox{\qquad$\Box$\hss}}

\newcommand{\tfrac}[2]{{\textstyle\frac{#1}{#2}}}

\makeatletter
\@addtoreset{equation}{section}
\makeatother
\newcommand{\eqref}[1]{{\rm (\ref{#1})}}

\newcommand{\sumcic}{\mathop{\mbox{\large {$\mathfrak S$}\vrule width 0pt
depth 2pt}}}

\newcommand{\End}{\mathop{\mbox{\rm End}}}
\newcommand{\Id}{\mathop{\mbox{\sl Id}}}
\newcommand{\SP}{\mathop{\mbox{\sl Sp}}}
\newcommand{\Un}{\mathop{\mbox{\sl U}}}
\newcommand{\un}{\mathop{\mbox{$\mathfrak u$}}}

\newcommand{\lcf}{\lbrack\!\lbrack}
\newcommand{\rcf}{\rbrack\!\rbrack}
\newcommand{\real}[1]{\lcf #1 \rcf}
\newcommand{\coderiv}{{\mkern 1mu d^*\mkern-1mu}}

\newcommand{\ncount}{167}
\newcommand{\ecount}{144}
\newcommand{\fcount}{7}
\newcommand{\nhkcount}{276}
\newcommand{\ehkcount}{316}

\newcommand{\nlchkcount}{44}
\newcommand{\elchkcount}{44}

\begin{document}
\title{Almost Hermitian Structures and Quaternionic Geometries}
\author{Francisco Mart\'\i n Cabrera and Andrew Swann}
\date{}
\markboth{\small {\sc Francisco Mart\'\i n Cabrera and Andrew Swann } }
         {\small \sc Almost Hermitian Structures and Quaternionic Geometries}
\baselineskip=5mm
\pagestyle{myheadings}

\maketitle

\begin{abstract}
  Gray \& Hervella gave a classification of almost Hermitian structures
  $(g,I)$ into $16$ classes.  We systematically study the interaction
  between these classes when one has an almost hyper-Hermitian structure
  $(g,I,J,K)$.  In general dimension we find at most $\ncount$ different
  almost hyper-Hermitian structures.  In particular, we obtain a number of
  relations that give hyperK\"aher or locally conformal hyperK\"ahler
  structures, thus generalising a result of Hitchin.  We also study the
  types of almost quaternion-Hermitian geometries that arise and tabulate
  the results.
\end{abstract}



\newpage
\tableofcontents

\newpage
\section{Introduction}

In \cite{Gray-H:16} Gray \& Hervella gave a classification of almost
Hermitian manifolds in terms of the covariant derivatives of the K\"ahler
$2$-form.  This derivative has special symmetries and may naturally be
decomposed into four components lying in spaces they called $\mathcal
W_1,\dots,\mathcal W_4$.  This gives $2^4=16$ different classes of almost
Hermitian manifolds determined by which components are non-zero.  From
these one may easily read off properties such as integrability of the
almost complex structure or closure of the K\"ahler form,
cf.~Table~\ref{tab:aH-class}.

An almost quaternion-Hermitian manifold locally possesses three almost
complex structures defining almost Hermitian structures with respect to a
common metric and satisfying the multiplicative identities of the imaginary
quaternions.  An analogue of the Gray-Hervella classification may be
obtained for such structures by considering the covariant derivative of a
certain four-form.  In general dimensions, this leads to $2^6=64$ classes,
which were recently described in detail in~\cite{Cabrera:aqH}.

It is a natural question to ask how these two classifications interact.
Indeed various results in this direction are already known, the most
celebrated being Hitchin's proof~\cite{Hitchin:Riemann-surface} that for a
manifold to be hyperK\"ahler it is sufficient that the three K\"ahler
$2$-forms be closed.  One first observation is that the space of covariant
derivatives of three arbitrary two-forms is about twice as large as the
space of covariant derivatives of a quaternionic four-form.  One should
therefore expect to find more relations than simply those arising from the
fact that the third almost complex structure is the product of the first
two.  In this paper, we find systematically all such relations between
these covariant derivatives.  With this in place it is an easy matter to
read off various consequences in the style of Hitchin's result and to
obtain generalisations.

After recalling definitions and the relevant representation theory
in~\S\ref{sec:prelim}, the key technical points of the paper may be found
in \S\ref{sec:decompositions}: a first decomposition of the covariant
derivative \( \nabla\omega_I \) is given in Lemma~\ref{lem:naome}, and this
is refined in Propositions \ref{prop:3} and~\ref{prop:E}.  Conclusions
regarding almost Hermitian types are drawn in~\S\ref{sec:hH}, whereas the
consequences for the quaternionic geometry may be found in~\S\ref{sec:qH}.
The paper ends with tables summarising all the results and a short
discussion of numbers of cases and possible examples.

\paragraph{Acknowledgements}
Andrew Swann is a member of \textsc{Edge}, Research Training Network
\textsc{hprn-ct-\oldstylenums{2000}-\oldstylenums{00101}}, supported by The
European Human Potential Programme.  He is grateful to the Department of
Fundamental Mathematics at the University of La Laguna for kind hospitality
whilst working on this project.  Francisco Mart\'\i n Cabrera wishes to
thank the Department of Mathematics and Computer Science at the University
of Southern Denmark for hospitality during the writing of this paper.

\begin{table}[tp]
  \centering
  \begin{tabular}{cll}
    \hline\strut
    $ \mathcal C $&Name&Characteristic property\\
    \hline\strut
    $ \{0\} $&K\"ahler&$\nabla\omega_I=0 $\\
    $\mathcal W_2$&almost K\"ahler, symplectic&$d\omega_I=0$\\
    $\mathcal W_4$&locally conformal K\"ahler&$d\omega_I=
    \alpha\wedge\omega_I\ne0$\\
    $\mathcal W_1+\mathcal W_2$&$(1,2)$-symplectic&$(d\omega_I)^{1,2}=0$\\
    $\mathcal W_3+\mathcal
    W_4$&integrable, Hermitian&$d\Lambda^{1,0} \subset
    \Lambda^{2,0}+\Lambda^{1,1}$\\
    $\mathcal W_1+\mathcal W_2+\mathcal
    W_3$&semi-K\"ahler, co-symplectic&$\coderiv \omega_I=0$\\ 
    \hline
  \end{tabular}
  \caption{Important classes~$ \mathcal C $ in the Gray-Hervella
    classification of almost Hermitian manifolds (of dimension at
    least~$6$).} 
  \label{tab:aH-class}
\end{table}

\section{Preliminaries}
\label{sec:prelim}

A $4n$-dimensional manifold $M$ $(n > 1)$ is said to be {\it almost
quaternion-Hermitian}, if $M$ is equipped with a Riemannian metric
$\langle \cdot,\cdot \rangle$ and a rank-three subbundle $\mathcal G$ of the
endomorphism bundle $\End TM$, such that locally $\mathcal G$ has an {\it
adapted basis} $I, J, K$ satisfying $I^2 = J^2= -1$ and $K= IJ = - JI$, and
$\langle A X, A Y \rangle = \langle X,Y \rangle $, for all $X,Y \in T_x M$
and $A =I,J,K$.  This is equivalent to saying that $M$ has a reduction of
its structure group to $\SP(n)\SP(1)$.  An almost quaternion-Hermitian
manifold with a global adapted basis is called an \emph{almost
hyper-Hermitian} manifold.

There are three local K\"ahler-forms $\omega_A (X,Y) = \langle X , A Y
\rangle$, $A=I,J,K$.  From these one may define a global, non-degenerate
four-form~$\Omega$, the {\it fundamental form}, by the local formula
\begin{equation}
  \label{eq:quafor}
  \Omega = \omega_I \wedge \omega_I + \omega_J \wedge \omega_J +
  \omega_K \wedge \omega_K.
\end{equation}

\subsection{Covariant Derivative of one K\"ahler Form}
Let $\nabla$ denote the Levi-Civita connection.  As the metric is almost
Hermitian with respect to $A=I,J,K$, one obtains \cite{Gray-H:16}
\begin{equation}
  \label{eq:GHo}
  \nabla_X \omega_A ( A Y , A Z) = - \nabla_X \omega_A (Y,Z).
\end{equation}
On the other hand the relation $I=JK$, implies that the covariant
derivative of $\omega_I$ is determined by those of $\omega_J$ and
$\omega_K$.  This may be expressed symmetrically by
\cite{Fernandez-Moreiras:symmetry,Cabrera:aqH}
\begin{displaymath}
  \nabla_X \omega_I ( J Y , KZ) + \nabla_X \omega_J ( K Y , I Z) + \nabla_X
  \omega_K ( I Y , J Z) = 0,
\end{displaymath}
or directly by
\begin{displaymath}
  \nabla_X \omega_I ( Y , Z ) = \nabla_X \omega_J ( K Y , Z) - \nabla_X
  \omega_K ( Y , J Z) .
\end{displaymath}
The two other versions of this equation obtained by cyclically permuting
$I,J,K$ also hold, and in this paper we will use such results without
further comment.

The following conventions will be used in the sequel.  If $b$ is a
$(0,s)$-tensor, we write
\begin{displaymath}
  \begin{array}{c}
    A_{(i)}b(X_1, \dots, X_i, \dots , X_s) = - b(X_1, \dots , AX_i, \dots ,
    X_s),\\
    Ab(X_1,\dots,X_s) = (-1)^sb(AX_1,\dots,AX_s),
  \end{array}
\end{displaymath}
for $A=I,J,K$.  We also consider the natural extension of~$\langle\cdot
,\cdot \rangle$ to $(0,s)$-tensors given by
\begin{displaymath}
  \langle a, b \rangle = \sum_{i_1, \dots , i_s=1}^{4n} a (
  e_{i_1}, \dots , e_{i_s})  b (e_{i_1}, \dots , e_{i_s}),
\end{displaymath}
where $\{ e_1, \dots , e_{4n} \}$ an orthonormal basis for $T_x M$.  The
notation $\{ e_1, \dots , e_{4n} \}$ will also denote the corresponding
dual basis of one-forms.  In some situations we will write $g$ for the
Riemannian metric $\langle\cdot ,\cdot \rangle$.

Using these conventions, we may write our expression for the covariant
derivative of $\omega_I$ as
\begin{equation}
  \label{eq:naomi}
  \nabla \omega_I = - K_{(2)} \nabla \omega_J + J_{(2)} \nabla \omega_K.
\end{equation}

\subsection{Representation Theory}
Our key tool for refining the expression~\eqref{eq:naomi}
for~$\nabla\omega_I$ is the representation theory of $\SP(n)$, $\SP(1)$
and~$\Un(1)$.  We will follow the $E$--$H$ formalism used in
\cite{Salamon:Invent,Swann:symplectiques,Swann:MathAnn} to denote
irreducible $\SP(n)\SP(1)$-modules.  Thus, $E$ is the fundamental
representation of $\SP(n)$ on $\mathbb C^{2n} \cong \mathbb H^n$ via left
multiplication by quaternionic matrices, and $H$ is the representation of
$\SP(1)$ on $\mathbb C^2 \cong \mathbb H$ given by $q.\zeta = \zeta
\overline{q}$, for $q \in \SP(1)$ and $\zeta \in H$.  An
$\SP(n)\SP(1)$-structure on a manifold $M$ gives rise to local bundles $E$
and $H$ associated to these representation and identifies $TM
\otimes_{\mathbb R} \mathbb C \cong E \otimes_\mathbb C H$.

On $E$, there is a $\SP(n)$-invariant complex symplectic form $\omega_E$
and a Hermitian inner product given by $\langle x , y \rangle_\mathbb C =
\omega_E (x , \widetilde{y})$, for all $x,y \in E$ and being $\widetilde{y}
= y j$ ($y \to \widetilde{y}$ is a quaternionic structure map on $E=\mathbb
C^{2n}$ considered as right complex vector space).  The mapping $x \to
x^{\omega} = \omega_E ( \cdot, x)$ gives us an identification of $E$ with
its dual $E^*$.  If $\{ u_1, \dots , u_n, \widetilde{u}_1, \dots ,
\widetilde{u}_n \}$ is a complex orthonormal basis for $E$, then
\begin{displaymath}
  \omega_E = u^{\omega}_i \wedge \widetilde{u}^{\omega}_i = u^{\omega}_i
  \widetilde{u}^{\omega}_i - \widetilde{u}^{\omega}_i u^{\omega}_i,
\end{displaymath}
where we have used the summation convention and omitted tensor product
signs.  These conventions will be used throughout the paper.

The irreducible representations of $\SP(1)$ are the symmetric powers $S^k H
\cong \mathbb C^{k+1}$.  An irreducible representation of $\SP(n)$ is
determined by its dominant weight $(\lambda_1,\dots,$ $ \lambda_n)$, where
$\lambda_i$ are integers with $\lambda_1 \geq \lambda_2 \geq \dots \geq
\lambda_n \geq 0$. This representation will be denoted by $V^{ (\lambda_1,
\dots , \lambda_r)}$, where $r$ is the largest integer such that $\lambda_r
> 0$.  Familiar notation is used for some of these modules when possible.
For instance, $V^{(k)} = S^k E$, the $k$th symmetric power of~$E$, and
$V^{(1, \dots , 1)} = \Lambda_0^r E$, where there are $r$ ones in the
exponent and $\Lambda_0^r E$ is the $\SP(n)$-invariant complement to $
\omega_E \Lambda^{r-2} E$ in $\Lambda^r E$.  Also $K$ will denote the
module $V^{(21)}$, which arises in the decomposition
\begin{displaymath}
  E \otimes \Lambda_0^2 E \cong \Lambda_0^3 E + K + E,
\end{displaymath}
where $+$ denotes direct sum.

When we fix a choice of local almost complex structure~$I$, we get a
reduction of the structure group to $\SP(n)\Un(1)\subset\Un(2n)$.  We write
$L$ for the standard representation of $\Un(1)$ on $\mathbb C$, and
$\Lambda^{1,0}$ for the representation of $\Un(2n)$ on $\mathbb C^{2n}$.
The latter has dual representation $\Lambda^{0,1}$, and arbitrary
irreducible representations lie in some tensor product
$\Lambda^{p,q}=\Lambda^p\Lambda^{1,0}\otimes\Lambda^q\Lambda^{0,1}$ and
will be labelled by the minimal pair~$(p,q)$.  In particular, $U^{3,0}$ is
the irreducible module in the decomposition
\begin{displaymath}
  \Lambda^{1,0}\otimes\Lambda^{2,0} = \Lambda^{3,0} + U^{3,0}.
\end{displaymath}
The space $\Lambda_0^{p,q}$ is the orthogonal complement of
$\omega_I\Lambda^{p-1,q-1}$ in $\Lambda^{p,q}$.  When we need to regard
these as representations over the real numbers we will use the notation
$\real{\Lambda_0^{p,q}}$, etc.  This satisfies
$\real{\Lambda_0^{p,q}}\otimes\mathbb C = \Lambda_0^{p,q}+\Lambda_0^{p,q}$.
Note that $TM=\real{\Lambda^{1,0}}$.  Elements of $\Lambda^{p,q}$ are said
to have type~$(p,q)$; elements of $\real{\Lambda^{p,q}}$ for $p\ne q$ will
be said to have type $\{p,q\}$.  Modules and types will be labelled by the
almost complex structure $I$ when it is necessary to avoid confusion.

\section{Quaternionic Decompositions}
\label{sec:decompositions}
Let us reconsider the covariant derivative of a single K\"ahler form
$\omega_I$, this time from the point of view of representation theory.
Equation~\eqref{eq:GHo} shows that
\begin{displaymath}
  \begin{array}{rccccccc}
  \nabla\omega_I \in \real{\Lambda^{1,0}}\otimes \un(2n)^\bot
  = &\real{\Lambda^{3,0}} &+ &\real{U^{3,0}} &+ &\real{\Lambda^{2,1}_0}
  &+ &\real{\Lambda^{1,0}},\\
  =&\mathcal W_1&+&\mathcal W_2&+&\mathcal W_3&+&\mathcal W_4,
\end{array}
\end{displaymath}
where $\un(2n)^\bot = \real{\Lambda^{2,0}}$ is the orthogonal complement of
$\un(2n)=\Lambda^{1,1}$ in $\Lambda^2T^*M$.  This is Gray \& Hervella's
decomposition~\cite{Gray-H:16} (cf.~\cite{Falcitelli-Farinola-Salamon}).

In the present context, we have also the action of $\SP(n)\Un(1)$ which is
a subgroup of $\Un(2n)$.  To obtain descriptions of the modules $\mathcal
W_i$ as representations of this subgroup, note that
\begin{displaymath}
  TM\otimes\mathbb C = EH = E(L+\overline L) = EL + E\overline L=
  \Lambda^{1,0} + \Lambda^{0,1}.
\end{displaymath}
From this we find
\begin{eqnarray}
  \label{eq:Lambda2}
    \Lambda^2T^*M\otimes\mathbb C &=& S^2E + S^2H + \Lambda^2_0ES^2H\\
    &=& S^2E + \mathbb C\omega_I + L^2 + \overline L^2 +
    \Lambda^2_0E(L^2+\mathbb C+\overline L^2) \nonumber
\end{eqnarray}
and hence $ \un(2n)\otimes\mathbb C = S^2E + \mathbb C\omega_I +
\Lambda^2_0E$, $\Lambda^{2,0} = (\Lambda^2_0E+\mathbb C)L^2$.  Taking the
tensor products of these representations we get
\begin{equation}
  \label{eq:dw} 
  \begin{array}{c}
  \mathcal W_1 \otimes \mathbb C =  (\Lambda_0^3 E + E)(L^3 +
  \overline{L}^3 ),\quad
  \mathcal W_2 \otimes \mathbb C  =   (K + E)(L^3  + \overline{L}^3 ),
  \\
  \mathcal W_3 \otimes \mathbb C =   (\Lambda_0^3 E + K  + E)(L +
  \overline{L}),\quad
  \mathcal W_4 \otimes \mathbb C  =   E(L + \overline{L}).
\end{array}
\end{equation}

We now set about finding the corresponding components of $\nabla\omega_I$
explicitly.  First note that $S^2H\subset \Lambda^2T^*M\otimes\mathbb C$
decomposes orthogonally both as $\mathbb C\omega_I + L^2 + \overline L^2$
and as $\mathbb C\omega_I+\mathbb C\omega_J+\mathbb C\omega_K$.  Thus
$L^2+\overline L^2$ is the direct sum of the tensors in $S^2H$ which are of
type $(1,1)_J$ or of type $(1,1)_K$.  Let us write $(\Lambda^2_0E)_I$ for
the module $\Lambda^2_0E$ in~\eqref{eq:Lambda2}; this consists of the
tensors in $\Lambda^2_0ES^2H$ which are of type $(1,1)_I$.  Our discussion
of $S^2H$ now shows $\Lambda^2_0E(L^2+\overline
L^2)=(\Lambda^2_0E)_J+(\Lambda^2_0E)_K$.  Thus
\begin{displaymath}
  \nabla_X\omega_I \in \un(n)^\bot = \mathbb R\omega_J + \mathbb R\omega_K
  + \real{(\Lambda^2_0E)_J} + \real{(\Lambda^2_0E)_K},
\end{displaymath}
with the last decomposition being orthogonal.

In order to label these components of $\nabla_X\omega_I$, consider the one
forms $\lambda_I$, $\lambda_J$, $\lambda_K$, defined by
\begin{equation}
  \label{eq:lambda-I}
  \lambda_I(X) := \tfrac1{4n} \langle \nabla_X \omega_J , \omega_K \rangle =
  - \tfrac1{4n} \langle \nabla_X \omega_K , \omega_J \rangle.
\end{equation}
For any two-form $\beta$ of type $(1,1)_I$, $I_{(1)}\beta$~is a symmetric
two-tensor.  Now
\begin{displaymath}
  S^2T^*M\otimes\mathbb C = S^2ES^2H + \Lambda^2_0E + \mathbb Cg.
\end{displaymath}
Thus for $\beta\in(\Lambda^2_0E)_I$, we have that
$I_{(1)}\beta\in\Lambda^2_0E$ and in particular $I_{(1)}\beta$ is of type
$(1,1)$ for $I$, $J$ \emph{and} $K$.  We therefore introduce the tensors
$\alpha_I$, $\alpha_J$, $\alpha_K$ in $T^*M\otimes\Lambda^2_0E\subset
T^*M\otimes S^2T^*M$ given by
\begin{eqnarray}
  \label{eq:alpha-I}
  \alpha_I &:=& - \lambda_I \otimes g + \tfrac12 ( J_{(2)} -J_{(3)})
  \nabla \omega_K \\
  &=& - \lambda_I \otimes g + \tfrac12 ( K_{(3)} -K_{(2)}) \nabla
  \omega_J.
  \nonumber
\end{eqnarray}
Note that they satisfy $\langle \alpha_A (X, \cdot , \cdot ),g\rangle = 0$
and
\begin{displaymath}
  \alpha_A=I_{(2)}I_{(3)} \alpha_A = J_{(2)}J_{(3)} \alpha_A =
  K_{(2)}K_{(3)} \alpha_A,
\end{displaymath}
for $A = I,J,K$.  These forms and tensors will play a significant r\^ole in
the present paper.

From equation~\eqref{eq:naomi} we have:

\begin{lemma}
  \label{lem:naome}
  Using the definitions \eqref{eq:lambda-I} and~\eqref{eq:alpha-I}, the
  covariant derivative of $\omega_I$ is given by
  \begin{displaymath}
    \nabla \omega_I = \lambda_K \otimes \omega_J - \lambda_J \otimes
    \omega_K + J_{(2)} \alpha_K - K_{(2)} \alpha_J.\qedhere
  \end{displaymath}
\end{lemma}

Since $\alpha_I, \alpha_J, \alpha_K \in EH \otimes \Lambda_0^2 E = (
\Lambda_0^3 E + K + E) H $, we may decompose $\alpha_I$ into three components
\begin{displaymath}
  \alpha_I = \alpha_I^{(3)} + \alpha_I^{(K)} + \alpha_I^{(E)} \in
  \Lambda_0^3 EH + KH + EH.
\end{displaymath}
(When $ \dim M=8 $, the module $ \Lambda^3_0E $ is trivial and the
corresponding component~$ \alpha_I^{(3)} $ is not present.)

Similar notation will be used for the $\SP(n)\Un(1)$-components in the
decomposition of~$\nabla\omega_I$.  Thus, for example,
$(\nabla\omega_I)_{\mathcal W_1}=(\nabla\omega_I)_{\mathcal
W_1}^{(3)}+(\nabla\omega_I)_{\mathcal W_1}^{(E)}$.  Controlling the
$\Lambda^3_0E$ and $K$ components of $\nabla\omega_I$ is relatively
straightforward as we shall now see.

The $\mathcal W_1 + \mathcal W_2$-part of $\nabla \omega_I$ consists of the
components that are of type $\{3,0\}_I$.  As $\nabla\omega_I$ is already of
type $\{2,0\}_I$ in the last two indices, one sees that
$2(\nabla\omega_I)_{\mathcal W_1 + \mathcal W_2} = (1- I_{(1)} I_{(2)})
\nabla \omega_I $.  Therefore, using Lemma~\ref{lem:naome}, we have
\begin{eqnarray}
  \label{eq:w1w2}
  \left( \nabla \omega_I \right)_{ \mathcal W_1 + \mathcal W_2 }
  & = &
  \tfrac12 \bigl\{ K ( J \lambda_J - K \lambda_K) \otimes \omega_J + J (
    J \lambda_J - K \lambda_K) \otimes \omega_K  \\
  && \qquad +(J_{(1)}  K_{(2)}   + K_{(1)}  J_{(2)}  ) (J_{(1)}
    \alpha_J - K_{(1)} \alpha_K ) \bigr\}. \nonumber
\end{eqnarray}
Similarly, the $\mathcal W_3 + \mathcal W_4$-component of $\nabla \omega_I$
is given by $2 ( \nabla \omega_I )_{\mathcal W_3 + \mathcal W_4} = (1+
I_{(2)} I_{(3)}) \nabla \omega_I $.  Thus, we have
\begin{eqnarray}
  \label{eq:w3w4}
  \left( \nabla \omega_I \right)_{ \mathcal W_3 + \mathcal W_4 } & = &
  \tfrac12 \bigl\{ - K ( J \lambda_J + K \lambda_K) \otimes \omega_J +
  J ( J \lambda_J + K \lambda_K) \otimes \omega_K \\
  && \qquad + (J_{(1)} K_{(2)} - K_{(1)} J_{(2)} ) ( J_{(1)} \alpha_J +
  K_{(1)} \alpha_K ) \bigr\}. \nonumber
\end{eqnarray}

Now for $\mathcal W_1 + \mathcal W_2$, the module $\Lambda^3_0E$ only
occurs in $\mathcal W_1$.  Similarly, in $\mathcal W_3 + \mathcal W_4$ this
module lies solely in $\mathcal W_3$.  Arguing in a similar way for the
$K$-components, we are lead to the following result.

\begin{proposition}
  \label{prop:3}
  Let $\alpha_I$ be as in~\eqref{eq:alpha-I}.  Then
  
  {\rm (a)} $(\nabla \omega_I)^{(3)}_{\mathcal W_1}$ is uniquely determined
  by $J_{(1)}\alpha_J^{(3)}-K_{(1)}\alpha_K^{(3)}$,
  
  {\rm (b)} $(\nabla \omega_I)^{(3)}_{\mathcal W_3}$
  by $J_{(1)}\alpha_J^{(3)}+K_{(1)}\alpha_K^{(3)}$,
  
  {\rm (c)} $(\nabla \omega_I)^{(K)}_{\mathcal W_2}$
  by $ J_{(1)}\alpha_J^{(K)}-K_{(1)}\alpha_K^{(K)}$ and
    
  {\rm (d)} $(\nabla \omega_I)^{(K)}_{\mathcal W_3}$ by $
  J_{(1)}\alpha_J^{(K)}+K_{(1)}\alpha_K^{(K)}$.
\end{proposition}

\begin{proof}
  The $\Lambda^3_0E$-parts of $(\nabla \omega_I)_{\mathcal W_1}$ and
  $(\nabla\omega_I)_{\mathcal W_3}$ are given by
  the corresponding parts of equations \eqref{eq:w1w2} and~\eqref{eq:w3w4}, i.e.,  \begin{eqnarray*}
    \left( \nabla \omega_I \right)_{\mathcal W_1}^{(3)} &=&
    \tfrac12(J_{(1)} K_{(2)} + K_{(1)} J_{(2)}) (J_{(1)} \alpha_J^{(3)}
    - K_{(1)} \alpha_K^{(3)}),\\
    \left( \nabla \omega_I \right)_{\mathcal W_3}^{(3)} &=& \tfrac12(J_{(1)}
    K_{(2)} - K_{(1)} J_{(2)}  ) ( J_{(1)} \alpha_J^{(3)} + K_{(1)}
    \alpha_K^{(3)} ),
  \end{eqnarray*}
  and analogous formul\ae{} for the $K$-components.  The Proposition now
  follows from the following Lemma. \cqd
\end{proof}

\begin{lemma}
  \label{jjaekka}
  For $\varepsilon=\pm1$, $(J_{(1)} K_{(2)} +\varepsilon K_{(1)} J_{(2)} )
  ( J_{(1)} \alpha_J - \varepsilon K_{(1)} \alpha_K ) = 0$ if and only if
  $J_{(1)} \alpha_J = \varepsilon K_{(1)} \alpha_K$.
\end{lemma}

\begin{proof}
  Let us just give the proof for $\varepsilon=+1$.  Clearly we only need to
  consider the forward implication.  Since $(J_{(1)} K_{(2)} + \varepsilon
  K_{(1)} J_{(2)} )^2 = 2 ( 1 + I_{(1)} I_{(2)})$, then the first equation
  of the lemma gives
  \begin{displaymath}
    J_{(1)} \alpha_J - K_{(1)} \alpha_K +  K_{(1)} I_{(2)} \alpha_J -
    J_{(1)} I_{(2)} \alpha_K = 0.
  \end{displaymath}
  But $ J_{(1)} \alpha_J - K_{(1)} \alpha_K \in T^* \otimes S^2 T^*$ and
  $K_{(1)} I_{(2)} \alpha_J - J_{(1)} I_{(2)} \alpha_K \in T^* \otimes
  \Lambda^2 T^*$.  Hence both of these pairs must vanish and in particular
  $J_{(1)} \alpha_J - K_{(1)} \alpha_K =0$ as required. \cqd
\end{proof}

The situation for the $E$-components of $\nabla\omega_I$ is a little more
complicated.  The presence of the $EH$-component in the decomposition of
$\alpha_I$ means that we can define one-forms $\eta_I$ from $\alpha_I$.  We
set
\begin{equation}
  \label{eq:eta-I}
  \eta_I (X) =  \alpha_I(e_i,e_i,X) = \langle \alpha_I(\cdot,\cdot,X), g
  \rangle.
\end{equation}
From a one-form $\eta$ we may produce a tensor~$\alpha(\eta)$ in $EH\otimes
\Lambda^2_0E$ by
\begin{displaymath}
  \alpha_I(\eta) = e_i\otimes(\eta\vee e_i + I\eta\vee Ie_i
  + J\eta\vee Je_i + K\eta\vee Ke_i) - \tfrac1n\eta\otimes g,
\end{displaymath}
where $a\vee b = \frac12(a\otimes b+b\otimes a)$.  As the map
$\eta\mapsto\alpha_I(\eta)$ is equivariant for the $\SP(n)\SP(1)$-action,
$\alpha_I(\eta)$ lies in $EH\subset EH\otimes\Lambda^2_0E$.  Computing the
corresponding one-form $\eta_I$ from $\alpha_I(\eta)$ via \eqref{eq:eta-I}
one gets $(2n+1)(n-1)\eta/n$.  Thus in general
\begin{equation}
  \label{eq:Ealpha}
  \alpha_I^{(E)} = \tfrac1{(2n+1)(n-1)}\bigl\{4n\, e_i\otimes(\eta_I\vee
  e_i)^{\mathbb H} - \eta_I\otimes
  g\bigr\},
\end{equation}
where $ 4(a\vee b)^{\mathbb H}= a\vee b + Ia\vee Ib + Ja\vee Jb + Ka\vee
Kb $.

\begin{proposition}
  \label{prop:E}
  Let $\lambda_I$ and $\eta_I$ be as in \eqref{eq:lambda-I}
  and~\eqref{eq:eta-I}.  Then
  
  {\rm (a)} $(\nabla \omega_I)^{(E)}_{\mathcal W_1}$ and $(\nabla
  \omega_I)^{(E)}_{\mathcal W_2}$ are uniquely determined by independent
  linear combinations of $J\lambda_J-K\lambda_K$ and $J\eta_J-K\eta_K$;
  
  {\rm (b)} $(\nabla \omega_I)^{(E)}_{\mathcal W_3}$ and $(\nabla
  \omega_I)^{(E)}_{\mathcal W_4}$ are uniquely determined by independent
  linear combinations of $J\lambda_J+K\lambda_K$ and $J\eta_J+K\eta_K$.
\end{proposition}

\begin{proof}
  For~(a), we compute the $E$-components of $ \left( \nabla \omega_I
  \right)_{ \mathcal W_1 + \mathcal W_2 }$ in the decompositions
  \eqref{eq:dw}.  For simplicity, write
  \begin{displaymath}
    \lambda^-_I :=  J \lambda_J - K \lambda_K, \qquad \eta^-_I :=  J \eta_J -
    K \eta_K.
  \end{displaymath}
  The $E$-component of $ \left( \nabla \omega_I \right)_{\mathcal W_1 +
  \mathcal W_2}$ is obtained from \eqref{eq:w1w2}, by replacing $\alpha_J$
  and $\alpha_K$ by $\alpha_J^{(E)}$ and $\alpha_K^{(E)}$, respectively.
  Taking \eqref{eq:Ealpha} into account, we have
  \begin{eqnarray*}
    2\left( \nabla \omega_I \right)_{ \mathcal W_1 + \mathcal W_2 }^{(E)}
    & = & K \left( \lambda_I^- + k_{12}\eta_I^-
    \right) \otimes \omega_J
    +  J \left( \lambda_I^- + k_{12} \eta_I^-
    \right) \otimes \omega_K \\
    && \qquad
    - n k_{12} \left( e_i \otimes Je_i \wedge K \eta_I^- + e_i \otimes Ke_i
    \wedge J \eta_I^- \right), 
  \end{eqnarray*}
  where $k_{12}= -1/(2n+1)(n-1)$.
  
  The $E$-part of the $\mathcal W_1$-component of $\nabla \omega_I $ is
  obtained by alternating $ \left( \nabla \omega_I \right)_{ \mathcal W_1 +
  \mathcal W_2 }^{(E)} $, i.e.,
  \begin{displaymath}
    6\left( \nabla \omega_I \right)_{ \mathcal W_1  }^{(E)} =
    K \left( \lambda_I^- + k_1\eta_I^- \right) \wedge \omega_J   +  J \left(
      \lambda_I^- + k_1\eta_I^- \right) \wedge \omega_K,
  \end{displaymath}
  where $k_1=-1/(n-1)$.  The difference $\left( \nabla \omega_I \right)_{
  \mathcal W_1 + \mathcal W_2 }^{(E)} - \left( \nabla \omega_I \right)_{
  \mathcal W_1 }^{(E)}$ is the $E$-part of the $\mathcal W_2$-component of
  $ \nabla \omega_I $, i.e.,
  \begin{eqnarray*}
    6\left( \nabla \omega_I \right)_{  \mathcal W_2 }^{(E)} & = & 2
    K \left( \lambda_I^- + k_2 \eta_I^- \right) \otimes \omega_J + 2 J \left(
      \lambda_I^- + k_2\eta_I^- \right) \otimes \omega_K \\
    && \qquad +   e_i \otimes Je_i \wedge K \left( \lambda_I^-  +
      k_2 \eta_I^- \right) + e_i \otimes Ke_i \wedge J \left( \lambda_I^- +
      k_2 \eta_I^- \right),
\end{eqnarray*}
where $k_2=1/(2n+1)$.  As $k_1\ne k_2$, this proves part~(a).

For part (b), we introduce
\begin{displaymath}
  \lambda_I^+  =  J \lambda_J + K \lambda_K, \qquad
  \eta_I^+  =  J \eta_J + K \eta_K.
\end{displaymath}
The $\mathcal W_4$-component of $\nabla \omega_I$ is given by~\cite{Gray-H:16}
\begin{displaymath}
  \left( \nabla \omega_I \right)_{\mathcal W_4} = \tfrac{1}{2(2n-1)} (- e_i
  \otimes e_i \wedge \coderiv  \omega_I + e_i \otimes Ie_i \wedge I
  \coderiv  \omega_I),
\end{displaymath}
where $\coderiv $ denotes the coderivative operator. Since from
Lemma~\ref{lem:naome} it follows that $I \coderiv \omega_I = \lambda_I^+ +
\eta_I^+$, then we have
\begin{displaymath}
   \left( \nabla \omega_I \right)_{\mathcal W_4} = \tfrac{1}{2(2n-1)} \left\{
     e_i \otimes e_i \wedge I (\lambda_I^+ + k_4\eta_I^+)  + e_i \otimes Ie_i
     \wedge (\lambda_I^+ + k_4\eta_I^+) \right\},
\end{displaymath}
where $ k_4=1 $.

The difference $\left( \nabla \omega_I \right)_{ \mathcal W_3 + \mathcal
W_4 }^{(E)} - \left( \nabla \omega_I \right)_{ \mathcal W_4 }^{(E)}$ is the
$E$-part of the $\mathcal W_3$-component of $ \nabla \omega_I $, i.e.,
\begin{eqnarray*}
  2\left( \nabla \omega_I \right)_{  \mathcal W_3 }^{(E)} & = & -
  K \left( \lambda_I^+ + k_3 \eta_I^+ \right) \otimes \omega_J +   J
  \left( \lambda_I^+ + k_3\eta_I^+ \right) \otimes \omega_K \\
  && - \tfrac{1}{(2n-1)} \left(e_i \otimes e_i \wedge I \left( \lambda_I^+
      +k_3\eta_I^+ \right) + e_i \otimes I e_i \wedge
  \left( \lambda_I^+ + k_3 \eta_I^+   \right)\right) ,
\end{eqnarray*}
where $k_3 = k_{12}=-1/(2n+1)(n-1)$ as before.  As $k_3\ne k_4$, we have
part~(b). \cqd
\end{proof}

\section{Almost Hyper-Hermitian Structures}
\label{sec:hH}
In the following theorem we show the way in which the class of $\omega_K$,
as an almost Hermitian structure, is conditioned by the respective classes
of $\omega_I$ and $\omega_J$.  Moreover, from the theorem one can also
deduce the list of possible triples of almost Hermitian types corresponding
to $\omega_I$, $\omega_J$ and $\omega_K$, respectively. Such a list is
contained in the tables given in~\S\ref{sec:tables}.  The theorem also
provides the essential rules for determining which triples of Gray-Hervella
types may occur.

\begin{theorem}
  \label{thm:clwww}
  Let $M$ be an almost hyper-Hermitian manifold.
  \begin{enumerate}
  \item[{\rm (i)}] If $\nabla \omega_I \in \mathcal W_3 + \mathcal W_4$
    and $\nabla \omega_J \in \mathcal C$, then $\nabla \omega_K \in \mathcal
    C$, where $\mathcal C$ means any Gray-Hervella class of almost
    Hermitian structures.
  \item[{\rm (ii)}] If $\nabla \omega_I, \nabla \omega_J \in \mathcal C +
    \mathcal W_3 + \mathcal W_4$, then $\nabla \omega_K \in \mathcal C+
    \mathcal W_3 + \mathcal W_4 $, for $\mathcal C= \mathcal W_1, \mathcal
    W_2$.
  \item[{\rm (iii)}] If $\nabla \omega_I,\nabla \omega_J \in \mathcal W_1 +
    \mathcal W_2$, then $\nabla \omega_K \in \mathcal W_3 + \mathcal W_4$.
  \item[{\rm (iv)}] If $\nabla \omega_I,\nabla \omega_J \in \mathcal W_1 +
    \mathcal W_2 + \mathcal W_4$ and $\nabla \omega_K \in \mathcal C + \mathcal W_3
    + \mathcal W_4$, for $\mathcal C = \mathcal W_1, \mathcal W_2$, then $\nabla
    \omega_K \in \mathcal W_3 + \mathcal W_4$.
  \end{enumerate}
  Moreover, if $M$ is eight-dimensional, we also have
  \begin{enumerate}
  \item[{\rm (v)}] If $\nabla \omega_I,\nabla \omega_J \in \mathcal W_1 +
    \mathcal W_2 + \mathcal W_3$ and $\nabla \omega_K \in \mathcal W_1 +
    \mathcal W_3 + \mathcal W_4$, then $\nabla \omega_K \in \mathcal W_3 +
    \mathcal W_4$.
  \end{enumerate}
\end{theorem}

\begin{proof}
  Follows directly from Propositions~\ref{prop:3} and~\ref{prop:E}.  \cqd
\end{proof}

We refer the reader to Table~\ref{tab:aH-class} for interpretations of some
of these classes.  In particular, parts~(i) and~(ii) each contain the
statement that if $I$ and $J$ are integrable then $K=IJ$ is too, as first
shown by Obata~\cite{Obata:connection}.

An important class of almost hyper-Hermitian manifolds are those in which
all three K\"ahler forms are parallel. These are hyperK\"ahler manifolds
and their metrics are Ricci-flat. The following result is a consequence of
the previous Theorem and shows some of the possible conditions which imply
that a manifold is hyperK\"ahler.

\begin{theorem}
  \label{thm:conseq}
  Let $M$ be an almost hyper-Hermitian manifold.  If one of the following
  conditions holds, then $M$~is hyperK\"ahler:
  \begin{enumerate}
  \item[{\rm (i)}] $\nabla \omega_I, \nabla \omega_J, \nabla \omega_K \in
    \mathcal W_1 + \mathcal W_2$,
  \item[{\rm (ii)}] $\nabla \omega_I \in \mathcal W_1$ and $\nabla \omega_J
    \in \mathcal W_2$,
  \item[{\rm (iii)}] $\nabla \omega_I =0$ and $\nabla \omega_J \in \mathcal
    W_1+\mathcal W_2$,
  \item[{\rm (iv)}] $\nabla \omega_I \in \mathcal C$ and $\nabla \omega_J \in
    \mathcal W_4$, $\mathcal C = \mathcal W_1, \mathcal W_2, \mathcal W_3$,
  \item[{\rm (v)}] $\nabla \omega_I=0$ and $\nabla \omega_J \in \mathcal C +
    \mathcal W_4$, $\mathcal C = \mathcal W_1, \mathcal W_2, \mathcal W_3$,
  \item[{\rm (vi)}] $\nabla \omega_I=0$, $\nabla \omega_J \in \mathcal W_1 +
    \mathcal W_2 + \mathcal W_4$ and $\nabla \omega_K \in \mathcal C + \mathcal
    W_3$, $\mathcal C= \mathcal W_1 + \mathcal W_2, \mathcal W_1 + \mathcal W_4,
    \mathcal W_2 + \mathcal W_4$, or
  \item[{\rm (vii)}] $\nabla \omega_I=0$, $\nabla \omega_J \in \mathcal W_1
    + \mathcal W_3 + \mathcal W_4$ and $\nabla \omega_K \in \mathcal W_2 +
    \mathcal W_3 + \mathcal W_4$,
  \end{enumerate}
  
  Moreover, if $M$ is eight-dimensional, each one of the following
  conditions also implies that $M$ is hyperK\"ahler:
  \begin{enumerate}
  \item[{\rm (viii)}] $\nabla \omega_I\in \mathcal W_1$ and $\nabla \omega_J
    \in \mathcal W_3$,
  \item[{\rm (ix)}] $\nabla \omega_I=0$ and $\nabla \omega_J \in \mathcal W_1
    + \mathcal W_3$, or
  \item[{\rm (x)}] $\nabla \omega_I=0$, $\nabla \omega_J \in \mathcal W_1 +
    \mathcal W_2+ \mathcal W_3$ and $\nabla \omega_K \in \mathcal W_1 + \mathcal
    W_3+ \mathcal W_4$.\cqd
  \end{enumerate}
\end{theorem}

\begin{observation}
  Part~(i) of Theorem~\ref{thm:conseq} was already proved in
  \cite{Murakoshi-SY:integrability} and is a generalisation of Hitchin's
  result~\cite{Hitchin:Riemann-surface} that if $\omega_I$, $\omega_J$ and
  $\omega_K$ are closed, then the manifold is hyperk\"ahler.  Part~(iii)
  includes the statement that a K\"ahler manifold is automatically
  hyperK\"ahler as soon as one additional two-form is closed.
\end{observation}

The $\mathcal W_4$-part of the covariant derivative of an almost Hermitian
structure is linearly determined by its Lee form~\cite{Gray-H:16} defined,
in the present context, by $\theta_A = 1/(2n-1) \, A \coderiv \omega_A$,
for $A=I,J,K$.  Below it will be shown that if the structures determined by
$I$, $J$, $K$ are locally conformal K\"ahler, then they have a common Lee
form.  Thus, in such a case, we can say that the manifold is locally
conformal hyperK\"ahler. In general, with $\lambda_I$ and $\eta_I$ as
in~\eqref{eq:lambda-I} and~\eqref{eq:eta-I}, we have the following result:

\begin{lemma}
  Let $M$ be an almost hyper-Hermitian manifold. The three almost Hermitian
  structures have a common Lee form if and only if
  \begin{displaymath}
    I \lambda_I + I \eta_I = J \lambda_J + J \eta_J = K \lambda_K + K
    \eta_K. 
  \end{displaymath}
\end{lemma}

Note this happens when $\left( \nabla \omega_I \right)_{ \mathcal W_1 +
\mathcal W_2 }^{(E)} = \left( \nabla \omega_J \right)_{ \mathcal W_1 +
\mathcal W_2 }^{(E)} =0 $, which is the case if $\nabla \omega_I, \nabla
\omega_J \in \mathcal W_3+\mathcal W_4$.

\begin{proof}
  From Lemma~\ref{lem:naome} we have $ I \coderiv \omega_I = J \lambda_J +
  K \lambda_K + J \eta_J + K \eta_K $, and the result follows.  \cqd
\end{proof}

Combining this Lemma with Theorem~\ref{thm:clwww} we obtain:

\begin{theorem}
  \label{thm:conseq2}
  Let $M$ be an almost hyper-Hermitian manifold.  If one of the following
  conditions holds, then $M$~is locally conformal hyperK\"ahler:
  \begin{enumerate}
  \item[{\rm (i)}] $\nabla \omega_I \in \mathcal W_4$ and $\nabla \omega_J
    \in \mathcal W_i + \mathcal W_4$, $i=1,2,3$,
  \item[{\rm (ii)}] $\nabla \omega_I, \nabla \omega_J \in \mathcal W_i +
    \mathcal W_4$ and $\nabla \omega_K \in \mathcal W_1 + \mathcal W_2 +
    \mathcal W_4$, $i=1,2,3$,
  \item[{\rm (iii)}] $\nabla \omega_I \in \mathcal W_i + \mathcal W_4$,
    $\nabla \omega_J \in \mathcal W_3 + \mathcal W_4$ and $\nabla \omega_K
    \in \mathcal W_j +\mathcal W_3+ \mathcal W_4$, $i,j=1,2,3$ and $i \neq
    j$,
  \item[{\rm (iv)}] $\nabla \omega_I \in \mathcal W_1 + \mathcal W_4$ and
    $\nabla \omega_J \in \mathcal W_2 + \mathcal W_4$ and $\nabla \omega_K
    \in \mathcal W_i + \mathcal W_j + \mathcal W_4$, $i,j=1,2,3$ and $i
    \neq j$, or
  \item[{\rm (v)}] $\nabla \omega_I\in \mathcal W_4$, $\nabla \omega_J \in
    \mathcal C + \mathcal W_4$ and $\nabla \omega_J \in \mathcal D+
    \mathcal W_4$, $\mathcal C, \mathcal D= \mathcal W_1 + \mathcal W_2,
    \mathcal W_1 + \mathcal W_3, \mathcal W_2 + \mathcal W_3$ and $\mathcal
    C \neq \mathcal D$.\cqd
  \end{enumerate}
\end{theorem}

\begin{observation}
  Part~(i) of Theorem~\ref{thm:conseq2} is a generalisation of Obata's
  result~\cite{Obata:Hermitian} that a K\"ahler structure with an
  additional integrable complex structure is hyperK\"ahler.
\end{observation}

\section{Almost Quaternion-Hermitian Structures}
\label{sec:qH}
Up to this point we have concentrated on the types of the almost Hermitian
structures $(g,I)$, $(g,J)$ and $(g,K)$.  However, these may be regarded as
coming from an adapted basis for an almost quaternion-Hermitian structure,
and in dimension at least~$12$ such a structure has one of $64$ possible
types determined by the covariant derivative of the fundamental
$4$-form~$\Omega$~\eqref{eq:quafor}.  It is therefore interesting to find
what consequences the three almost Hermitian types have for the
quaternionic type.  The $64$ quaternionic classes come from the following
$\SP(n)\SP(1)$-decomposition.

\begin{proposition}[Swann~\cite{Swann:symplectiques}]
  The covariant derivative of the fundamental form $\Omega$ of an almost
  quaternion-Hermitian manifold $M$ of dimension at least $8$, has the
  property
  \begin{displaymath}
  \nabla \Omega \in T^* M \otimes \Lambda_0^2E S^2 H = ( \Lambda_0^3 E +
  K + E) ( S^3 H + H).\qedhere
  \end{displaymath}
\end{proposition}

\noindent
If the dimension of $M$ is at least~$12$, all the modules of the sum are
non-zero.  For an eight-dimensional manifold $M$, we have $\Lambda_0^3 E
S^3 H = \Lambda_0^3 E H = \{ 0 \}$.

If $\nabla\Omega=0$, $M$~is said to be quaternionic K\"ahler and the metric
$g$~is automatically Einstein (see for example \cite{Besse:Einstein}).  If
$\nabla\Omega\in EH$, then $M$~is locally conformal quaternionic
K\"ahler.  The case $\nabla\Omega\in (K+E)H$ is known as QKT geometry:
there is a second $\SP(n)\SP(1)$-connection on~$M$ with totally
skew-symmetric torsion, see for example~\cite{Ivanov:QKT}.  When
$\nabla\Omega\in (\Lambda^3_0E+K+E)H$ the underlying almost quaternionic
structure is integrable, i.e., there is a torsion-free connection
preserving the bundle spanned by $I$, $J$ and~$K$.

Using Lemma~\ref{lem:naome}, the covariant derivative $\nabla \Omega$ is
given by
\begin{equation}
  \label{eq:ornaomeg}
  \nabla \Omega = 2 \sumcic_{I,J,K} \bigl( J_{(2)} \alpha_K - K_{(2)}
    \alpha_J \bigr) \wedge \omega_I,
\end{equation}
where $\sumcic$ denotes the cyclic sum.  Note that the one-forms
$\lambda_I$, $\lambda_J$ and~$\lambda_K$ do not appear in this formula.  We
immediately conclude that the $\Lambda^3_0E(S^3H+H)$, $K(S^3H+H)$ and
$E(S^3H+H)$ of $\nabla\Omega$ are linearly determined by the corresponding
components of the $\alpha$'s.  To further divide these components we use
the following result.

Let ${\bf a} \colon \bigotimes^4 T^* M \to \Lambda^4 T^* M$ be the
alternation map.

\begin{proposition}[Cabrera~\cite{Cabrera:aqH}]
  The covariant derivative of $\Omega$ splits as
  \begin{displaymath}
      \nabla\Omega = (\nabla\Omega)_{S^3H} + (\nabla\Omega)_H
       \in  (\Lambda^3_0E+K+E)S^3H + (\Lambda^3_0E+K+E)H,
   \end{displaymath}
   with
   \begin{eqnarray}
     \label{eq:S3H}
     (\nabla\Omega)_{S^3H} &=& \tfrac16(4\nabla\Omega-\sumcic_{I,J,K} (\Id
     \otimes {\bf a}) I_{(1)}I_{(2)} \nabla \Omega),\\
     \label{eq:H}
     (\nabla\Omega)_H &=& \tfrac16(2\nabla\Omega+\sumcic_{I,J,K} (\Id
     \otimes {\bf a}) I_{(1)}I_{(2)} \nabla \Omega).\qedhere
   \end{eqnarray}
\end{proposition}

The classes involved in this last result are also determined by conditions
on $\alpha_I$, $\alpha_J$, $\alpha_K$.

\begin{proposition}
  Let $M$ an almost quaternion-Hermitian manifold and $U$ an open set where
  the adapted basis $I,J,K$ is defined. Then
  \begin{enumerate}
  \item[{\rm (i)}] $(\nabla\Omega)_H$ is linearly
    determined by
    $I_{(1)}\alpha_I+J_{(1)}\alpha_J+K_{(1)}\alpha_K$,
  \item[{\rm (ii)}] $(\nabla\Omega)_{S^3H}$ is linearly determined by
    $I_{(1)}\alpha_I-J_{(1)}\alpha_J$ and
    $J_{(1)}\alpha_J-K_{(1)}\alpha_K$.
  \end{enumerate}
\end{proposition}

\begin{proof}
  Using~\eqref{eq:ornaomeg} the tensors appearing in equations
  \eqref{eq:S3H} and~\eqref{eq:H} may be expressed as
  \begin{eqnarray*}
    \sumcic_{I,J,K} (\Id \otimes {\bf a}) I_{(1)}I_{(2)} \nabla
    \Omega   = 4 \sumcic_{I,J,K} && \bigl\{   J_{(1)} K_{(2)} (I_{(1)}
    \alpha_I + K_{(1)} \alpha_K) \\
    && \quad - K_{(1)}J_{(2)} ( I_{(1)} \alpha_I + J_{(1)} \alpha_J)
    \bigr\} \wedge \omega_I,
  \end{eqnarray*}
  \begin{displaymath}
    \nabla \Omega  =  2 \sumcic_{I,J,K} \bigl( - K_{(1)} J_{(2)} K_{(1)}
      \alpha_K + J_{(1)} K_{(2)}
      J_{(1)}\alpha_J \bigr) \wedge \omega_I.
  \end{displaymath}
  This gives
  \begin{displaymath}
    \tfrac32(\nabla\Omega)_H = \sumcic_{I,J,K}
    (J_{(1)}K_{(2)}-K_{(1)}J_{(2)})
    (I_{(1)}\alpha_I+J_{(1)}\alpha_J+K_{(1)}\alpha_K)
    \wedge \omega_I,
  \end{displaymath}
  \begin{eqnarray*}
    \tfrac32(\nabla\Omega)_{S^3H} = \sumcic_{I,J,K} && \bigl\{
    J_{(1)}K_{(2)}(2J_{(1)}\alpha_J-K_{(1)}\alpha_K-I_{(1)}\alpha_I)\\
    &&- K_{(1)}J_{(2)}(2K_{(1)}\alpha_K-I_{(1)}\alpha_I-J_{(1)}\alpha_J)
    \bigr\}\wedge\omega_I.
  \end{eqnarray*}
  Both of these expressions have the form
  \begin{displaymath}
    \Xi_\beta := \sumcic_{I,J,K}
    (K_{(2)}\beta_J-J_{(2)}\beta_K)\wedge\omega_I,
  \end{displaymath}
  with $ \beta_I \in T^*M\otimes \Lambda^2_0E \subset T^*M\otimes S^2T^*M $.

  Now, we compute
  \begin{displaymath}
    \sum_{i=1}^{4n} \Xi_\beta(X,Y,Z,e_i,Ie_i) =
    (4n+2)\bigl(\beta_J(X,KY,Z)-\beta_K(X,JY,Z)\bigr).
  \end{displaymath}
  However, $ I_{(2)}\beta_J(X, \cdot , \cdot) $ is an element of $
  \Lambda^2T^*M $ that is of type $ (1,1)_I $, $ \{2,0\}_J $ and~$
  \{2,0\}_K $.  So $ \Xi_\beta $ is uniquely determined by $ \beta_I $, $
  \beta_J $ and~$ \beta_K $ and the result follows.  \cqd
\end{proof}

Combining this Proposition with the remark after
equation~\eqref{eq:ornaomeg}, we find that for $ V=\Lambda^3_0E, K, E $ the
$ VH $-component of $ \nabla\Omega $ is uniquely determined by
$I_{(1)}\alpha^{(V)}_I+J_{(1)}\alpha^{(V)}_J+K_{(1)}\alpha^{(V)}_K$, and
the $ VS^3H $-component is uniquely determined by
$I_{(1)}\alpha^{(V)}_I-J_{(1)}\alpha^{(V)}_J$ and
$J_{(1)}\alpha^{(V)}_J-K_{(1)}\alpha^{(V)}_K$.  We can thus fully
determined the quaternionic type of the manifold from information about the
$ \alpha $'s.  As these are determined by the covariant derivatives
$\nabla\omega_I$, etc., we obtain the following relations between Hermitian
and quaternionic types (see also the tables in~\S\ref{sec:tables}).

\begin{theorem}
  Let $M$ an almost quaternion-Hermitian manifold.  On an open set~$U$
  where an adapted basis $I,J,K$ is defined, one has:
  \begin{enumerate}
  \item[{\rm (i)}] If $\nabla \omega_I \in \mathcal W_1 + \mathcal W_4$ and
    $\nabla \omega_J \in \mathcal W_2 + \mathcal W_4$, or $\nabla \omega_I,
    \nabla \omega_J, \nabla \omega_K \in \mathcal W_1+\mathcal W_2 +
    \mathcal W_4$, then $\nabla \Omega \in E (S^3 H + H)$.
  \item[{\rm (ii)}] If $\nabla \omega_I, \nabla \omega_J \in \mathcal
    W_1+\mathcal W_3 + \mathcal W_4$, then $\nabla \Omega \in ( \Lambda_0^3
    E + E) (S^3 H + H)+KH$.
  \item[{\rm (iii)}] If $\nabla \omega_I, \nabla \omega_J \in \mathcal
    W_2+\mathcal W_3 + \mathcal W_4$, then $\nabla \Omega \in ( K + E) (S^3
    H + H) + \Lambda_0^3 E H$.
  \item[{\rm (iv)}] If $\nabla \omega_I \in \mathcal W_1 + \mathcal W_4$
    and $\nabla\omega_J \in \mathcal W_1 + \mathcal C + \mathcal W_4$,
    $\mathcal C= \mathcal W_2,\mathcal W_3$, then $\nabla \Omega \in (
    \Lambda_0^3 E + E) (S^3 H + H)$.
  \item[{\rm (v)}] If $\nabla \omega_I \in \mathcal W_2 + \mathcal W_4$ and
    $\nabla\omega_J \in \mathcal W_2 + \mathcal C + \mathcal W_4$,
    $\mathcal C= \mathcal W_1,\mathcal W_3$, then $\nabla \Omega \in ( K +
    E) (S^3 H + H)$.
  \item[{\rm (vi)}] If $\nabla \omega_I, \nabla \omega_J \in \mathcal W_3 +
    \mathcal W_4$, then $\nabla \Omega \in ( \Lambda_0^3 E+ K + E) H$.
  \item[{\rm (vii)}] If $\nabla \omega_I \in \mathcal W_1$ and $\nabla
    \omega_J, \nabla \omega_K \in \mathcal W_1 + \mathcal W_3$, then
    $\nabla \Omega \in \Lambda_0^3 E(S^3 H+ H)$.
  \item[{\rm (viii)}] If $\nabla \omega_I \in \mathcal W_2$ and $\nabla
    \omega_J, \nabla \omega_K \in \mathcal W_2 + \mathcal W_3$, then
    $\nabla \Omega \in K(S^3 H+ H)$.
  \item[{\rm (ix)}] If $\nabla \omega_I,\nabla \omega_J, \nabla \omega_K
    \in \mathcal W_1 + \mathcal W_3$, then $\nabla \Omega \in \Lambda_0^3E
    (S^3 H+ H)+ (K+E)H$.
  \item[{\rm (x)}] If $\nabla \omega_I,\nabla \omega_J, \nabla \omega_K \in
    \mathcal W_2 + \mathcal W_3$, then $\nabla \Omega \in K(S^3 H+ H)+
    (\Lambda_0^3E+E)H$.
  \item[{\rm (xi)}] If one of the following conditions holds, then $\nabla
    \Omega \in (\Lambda_0^3E + K)(S^3 H+ H)$:
    \begin{enumerate}
    \item[{\rm (a)}] $\nabla \omega_I \in \mathcal W_1$, $\nabla \omega_J
      \in \mathcal W_2 + \mathcal W_3 $ and $\nabla \omega_K \in \mathcal
      W_1 + \mathcal W_2 + \mathcal W_3$,
    \item[{\rm (b)}] $\nabla \omega_I \in \mathcal W_2$, $\nabla \omega_J
      \in \mathcal W_1 + \mathcal W_3 $ and $\nabla \omega_K \in \mathcal
      W_1 + \mathcal W_2 + \mathcal W_3$, or
    \item[{\rm(c)}] $\nabla \omega_I \in \mathcal W_1+ \mathcal W_2$,
      $\nabla \omega_J \in \mathcal W_1 + \mathcal W_3 $ and $\nabla
      \omega_K \in \mathcal W_2 + \mathcal W_3$.
    \end{enumerate}
  \item[{\rm (xii)}] If $\nabla \omega_I \in \mathcal W_1 + \mathcal W_3$,
    $\nabla \omega_J \in \mathcal W_2 + \mathcal W_3$ and $\nabla \omega_K
    \in \mathcal W_1+ \mathcal W_2 + \mathcal W_3$, then $\nabla \Omega \in
    ( \Lambda_0^3 E+ K) (S^3 H + H) + EH$.
  \end{enumerate}

  Moreover, if $M$ is eight-dimensional, one also has
  \begin{enumerate}
  \item[{\rm (xiii)}] If $\nabla \omega_I \in \mathcal W_1 + \mathcal W_3$,
    $\nabla \omega_J \in \mathcal W_1 + \mathcal W_2 + \mathcal W_3$ and
    $\nabla \omega_K \in \mathcal W_2 + \mathcal W_4$, then $\nabla \Omega
    \in K (S^3 H + H) + E S^3 H$.
  \item[{\rm (xiv)}] If $\nabla \omega_I,\nabla \omega_J \in \mathcal W_1 +
    \mathcal W_3$ and $\nabla \omega_K \in \mathcal W_4$, then $\nabla
    \Omega \in E S^3 H$.\cqd
  \end{enumerate}
\end{theorem}

\section{Tables and Comments}
\label{sec:tables}

Tables \ref{tab:n1} and~\ref{tab:n2} show the full consequences of the
formul{\ae} derived in this paper.

Hermitian types are denoted by a hexadecimal number
$0$,\dots,$9$,$A$,\dots,$F$, where $W_i$ contributes $2^{i-1}$. So, for
example $B=1+2+8$ represents $\mathcal W_1+\mathcal W_2+\mathcal W_4$.

The rows of the table give the Hermitian type of $I$, the columns the type
of~$J$.  Due to symmetry we only need to show the cases where the Hermitian
type of~$J$ is greater than or equal to that of~$I$.

Each rectangle in the table contains up to $16$ entries corresponding to
the Hermitian types of~$K$ that are greater than or equal to that of~$J$.
These are arranged with the type of $K$ increasing in each column, so the
first column potentially begins with type~$0$, the next type~$4$, and then
type $8$ and finally type~$C$.

In each position in this rectangle there is one of two types of entry.
Three hexadecimal digits $abc$ indicate that the Hermitian types of $I$,
$J$ and $K$ reduce to $a$, $b$ and $c$ respectively.  Two bold digits
\textbf{PQ}, indicate that the Hermitian types do not reduce and specify
instead the quaternionic type of the manifold.  \textbf{P}~corresponds to
the~$S^3H$-part of $\nabla\Omega$ and \textbf{Q} to the $H$-part, with
$\Lambda^3_0E$~contributing~$4$, $K$~$2$ and $E$~$1$.  Thus \textbf{36}
indicates type $(K+E)S^3H+(\Lambda^3_0E+K)H$.  Note that the bottom right
entry in each rectangle corresponds to $\omega_K$ having type~$F=\mathcal
W_1+\mathcal W_2+\mathcal W_3+\mathcal W_4$; this is no restriction on
$\omega_K$ and so this entry tells what happens when $I$ and~$J$ have a
specified type.

These results are for general dimension $4n\geqslant12$.  In dimension~$8$,
several conclusions may be different.  These are indicated by italicising
the entry in Tables \ref{tab:n1} and~\ref{tab:n2}, provided the difference
does not simply arise because of the absence of $\Lambda^3_0E$ in the
quaternionic type.  What actually happens in dimension $8$ in these special
cases is then given in Table~\ref{tab:d8}.  This table lists the
$I$,$J$,$K$ type, its reduction and finally the quaternionic type.  The
entries different from the general case are again italicised.

One finds that in general dimension there are $\ncount$ different almost
hyper-Hermitian types, whilst in dimension~$8$ there are only $\ecount$.
In comparison, the number of potential triples of types is $\tfrac16
16.17.18=816$.  Of these $816$ cases, $\nhkcount$~are hyperK\"ahler and
$\nlchkcount$ are locally conformal hyperK\"ahler in general dimension.
For dimension $8$, one gets instead $\ehkcount$ and $\elchkcount$,
respectively.

For completeness Table~\ref{tab:d4} gives the situation for dimension~$4$.
In this case there are no $\alpha_I$ terms and the only Hermitian types are
$\{0\}$, $\mathcal W_2$, $\mathcal W_4$ and $\mathcal W_2+\mathcal W_4$.
We do not specify quaternionic types in this table as these are no longer
determined by~$\nabla\Omega$ (the four-form $\Omega$~is a constant times
the volume form, and so parallel).  We see that there are only $\fcount$
distinct almost hyper-Hermitian types in this case.

It is natural to ask whether examples of each of the $\ncount$ different
almost hyper-Hermitian types occur.  With so many cases this is clearly a
daunting task.  However, one special case that is of interest is when $I$,
$J$ and~$K$ have the same type.  In this situation one may check that if a
given component of $\nabla\omega_I$ vanishes then the same is true of
$\nabla\omega_A$, where $A=aI+bJ+cK$, with $a^2+b^2+c^2=1$ constant.  The
table shows that the only possible Gray-Hervella types are $\{0\}$,
$\mathcal W_3+\mathcal C$ (with $\mathcal C\subset\mathcal W_1+\mathcal
W_2+\mathcal W_4$), $\mathcal W_4$ or $\mathcal W_1+\mathcal W_2+\mathcal
W_4$.  With the exception of the last case, these may all be realised in
dimension~$12$ by considering homogeneous structures, and conformal changes
of such, on $(S^3)^4$, $T^3\times M^3$, with $M^3$ a three-dimensional Lie
group, either semi-simple, nilpotent or solvable.  However, the given
structures do not exhibit the full predicted almost quaternionic-Hermitian
types, and we have not yet found examples of the last case $\mathcal
W_1+\mathcal W_2+\mathcal W_4$, which will be quaternionic type
$E(S^3H+H)$.  We therefore reserve presentation of such examples to future
work.

\begin{table}[p]
  \centering
\newcommand{\tblfnt}{\fontsize{7}{8}\selectfont}%
\newcommand{\rowbox}[1]{\hbox to 18pt{\tblfnt \hfill\strut #1\hfill}}%
\newcommand{\colbox}[2][32pt]{\vbox to #1{\vfill #2}}%
\setlength{\tabcolsep}{1pt}%
\begin{tabular}{r|c|c|c|c|c|c|c|c|c|c|}
\setlength{\parindent}{0pt}
\hbox to 0pt{\hss$I\backslash$\hskip-3pt\raisebox{4pt}{$J$}}&F&E&D&C&B&A&9&8&7&6\\
\hline
\hbox{\colbox{\hbox{0 }\vfill}}
&\hbox{\colbox{\rowbox{\textbf{77}}}}
&\hbox{\colbox{\rowbox{\textbf{33}}\rowbox{0EE}}}
&\hbox{\colbox{\rowbox{\textbf{55}}\rowbox{000}\rowbox{0DD}}}
&\hbox{\colbox{\rowbox{000}\rowbox{000}\rowbox{000}\rowbox{000}}}
&\hbox{\colbox{\rowbox{\textbf{11}}}}%
\hbox{\colbox{\rowbox{000}\rowbox{000}\rowbox{000}\rowbox{0BB}}}
&\hbox{\colbox{\rowbox{000}\rowbox{000}}}%
\hbox{\colbox{\rowbox{000}\rowbox{000}\rowbox{000}\rowbox{000}}}
&\hbox{\colbox{\rowbox{000}\rowbox{000}\rowbox{000}}}%
\hbox{\colbox{\rowbox{000}\rowbox{000}\rowbox{000}\rowbox{000}}}
&\hbox{\colbox{\rowbox{000}\rowbox{000}\rowbox{000}\rowbox{000}}}%
\hbox{\colbox{\rowbox{000}\rowbox{000}\rowbox{000}\rowbox{000}}}
&\hbox{\colbox{\rowbox{\textbf{77}}}}%
\hbox{\colbox{\rowbox{000}\rowbox{000}\rowbox{000}\rowbox{000}}}%
\hbox{\colbox{\rowbox{000}\rowbox{\textit{055}}\rowbox{066}\rowbox{077}}}
&\hbox{\colbox{\rowbox{\textbf{22}}\rowbox{066}}}%
\hbox{\colbox{\rowbox{000}\rowbox{000}\rowbox{000}\rowbox{000}}}%
\hbox{\colbox{\rowbox{000}\rowbox{000}\rowbox{066}\rowbox{066}}}
\\
\hline
\hbox{\colbox{\hbox{1 }\vfill}}
&\hbox{\colbox{\rowbox{\textbf{77}}}}
&\hbox{\colbox{\rowbox{0EE}\rowbox{\textbf{77}}}}
&\hbox{\colbox{\rowbox{\textbf{55}}\rowbox{11C}\rowbox{1DD}}}
&\hbox{\colbox{\rowbox{000}\rowbox{1C1}\rowbox{000}\rowbox{1C1}}}
&\hbox{\colbox{\rowbox{0BB}}}%
\hbox{\colbox{\rowbox{11C}\rowbox{11C}\rowbox{11C}\rowbox{\textbf{55}}}}
&\hbox{\colbox{\rowbox{000}\rowbox{000}}}%
\hbox{\colbox{\rowbox{000}\rowbox{000}\rowbox{000}\rowbox{\textbf{11}}}}
&\hbox{\colbox{\rowbox{000}\rowbox{000}\rowbox{000}}}%
\hbox{\colbox{\rowbox{11C}\rowbox{11C}\rowbox{11C}\rowbox{11C}}}
&\hbox{\colbox{\rowbox{000}\rowbox{000}\rowbox{000}\rowbox{000}}}%
\hbox{\colbox{\rowbox{000}\rowbox{000}\rowbox{000}\rowbox{000}}}
&\hbox{\colbox{\rowbox{\textit{\textbf{77}}}}}%
\hbox{\colbox{\rowbox{000}\rowbox{\textit{141}}\rowbox{000}\rowbox{\textbf{55}}}}%
\hbox{\colbox{\rowbox{11C}\rowbox{\textit{15D}}\rowbox{\textbf{77}}\rowbox{\textbf{77}}}}
&\hbox{\colbox{\rowbox{066}\rowbox{\textit{\textbf{66}}}}}%
\hbox{\colbox{\rowbox{000}\rowbox{\textit{141}}\rowbox{000}\rowbox{\textit{141}}}}%
\hbox{\colbox{\rowbox{000}\rowbox{\textit{141}}\rowbox{066}\rowbox{\textbf{77}}}}
\\
\hline
\hbox{\colbox{\hbox{2 }\vfill}}
&\hbox{\colbox{\rowbox{\textbf{77}}}}
&\hbox{\colbox{\rowbox{\textbf{33}}\rowbox{2EE}}}
&\hbox{\colbox{\rowbox{0DD}\rowbox{2C2}\rowbox{\textbf{77}}}}
&\hbox{\colbox{\rowbox{000}\rowbox{000}\rowbox{2C2}\rowbox{2C2}}}
&\hbox{\colbox{\rowbox{0BB}}}%
\hbox{\colbox{\rowbox{22C}\rowbox{22C}\rowbox{22C}\rowbox{\textbf{33}}}}
&\hbox{\colbox{\rowbox{000}\rowbox{000}}}%
\hbox{\colbox{\rowbox{22C}\rowbox{22C}\rowbox{22C}\rowbox{22C}}}
&\hbox{\colbox{\rowbox{000}\rowbox{000}\rowbox{000}}}%
\hbox{\colbox{\rowbox{000}\rowbox{000}\rowbox{000}\rowbox{\textbf{11}}}}
&\hbox{\colbox{\rowbox{000}\rowbox{000}\rowbox{000}\rowbox{000}}}%
\hbox{\colbox{\rowbox{000}\rowbox{000}\rowbox{000}\rowbox{000}}}
&\hbox{\colbox{\rowbox{\textbf{77}}}}%
\hbox{\colbox{\rowbox{000}\rowbox{000}\rowbox{242}\rowbox{\textbf{33}}}}%
\hbox{\colbox{\rowbox{22C}\rowbox{\textit{\textbf{77}}}\rowbox{26E}\rowbox{\textbf{77}}}}
&\hbox{\colbox{\rowbox{\textbf{22}}\rowbox{266}}}%
\hbox{\colbox{\rowbox{000}\rowbox{000}\rowbox{242}\rowbox{242}}}%
\hbox{\colbox{\rowbox{22C}\rowbox{22C}\rowbox{\textbf{33}}\rowbox{26E}}}
\\
\hline
\hbox{\colbox{\hbox{3 }\vfill}}
&\hbox{\colbox{\rowbox{\textbf{77}}}}
&\hbox{\colbox{\rowbox{2EE}\rowbox{\textbf{77}}}}
&\hbox{\colbox{\rowbox{1DD}\rowbox{\textbf{77}}\rowbox{\textbf{77}}}}
&\hbox{\colbox{\rowbox{000}\rowbox{1C1}\rowbox{2C2}\rowbox{3C3}}}
&\hbox{\colbox{\rowbox{\textbf{11}}}}%
\hbox{\colbox{\rowbox{33C}\rowbox{33C}\rowbox{33C}\rowbox{\textbf{77}}}}
&\hbox{\colbox{\rowbox{000}\rowbox{383}}}%
\hbox{\colbox{\rowbox{22C}\rowbox{22C}\rowbox{22C}\rowbox{\textbf{33}}}}
&\hbox{\colbox{\rowbox{000}\rowbox{000}\rowbox{383}}}%
\hbox{\colbox{\rowbox{11C}\rowbox{11C}\rowbox{11C}\rowbox{\textbf{55}}}}
&\hbox{\colbox{\rowbox{000}\rowbox{000}\rowbox{000}\rowbox{383}}}%
\hbox{\colbox{\rowbox{000}\rowbox{000}\rowbox{000}\rowbox{383}}}
&\hbox{\colbox{\rowbox{\textbf{77}}}}%
\hbox{\colbox{\rowbox{338}\rowbox{\textit{\textbf{55}}}\rowbox{\textbf{33}}\rowbox{\textbf{77}}}}%
\hbox{\colbox{\rowbox{33C}\rowbox{\textit{\textbf{77}}}\rowbox{\textbf{77}}\rowbox{\textbf{77}}}}
&\hbox{\colbox{\rowbox{266}\rowbox{\textbf{77}}}}%
\hbox{\colbox{\rowbox{000}\rowbox{\textit{141}}\rowbox{242}\rowbox{343}}}%
\hbox{\colbox{\rowbox{22C}\rowbox{\textit{\textbf{77}}}\rowbox{26E}\rowbox{\textbf{77}}}}
\\
\hline
\hbox{\colbox{\hbox{4 }\vfill}}
&\hbox{\colbox{\rowbox{\textbf{77}}}}
&\hbox{\colbox{\rowbox{\textbf{37}}\rowbox{4EE}}}
&\hbox{\colbox{\rowbox{\textbf{57}}\rowbox{444}\rowbox{4DD}}}
&\hbox{\colbox{\rowbox{444}\rowbox{444}\rowbox{444}\rowbox{444}}}
&\hbox{\colbox{\rowbox{\textbf{77}}}}%
\hbox{\colbox{\rowbox{000}\rowbox{499}\rowbox{4AA}\rowbox{4BB}}}
&\hbox{\colbox{\rowbox{\textbf{33}}\rowbox{4AA}}}%
\hbox{\colbox{\rowbox{000}\rowbox{000}\rowbox{4AA}\rowbox{4AA}}}
&\hbox{\colbox{\rowbox{\textbf{55}}\rowbox{000}\rowbox{499}}}%
\hbox{\colbox{\rowbox{000}\rowbox{499}\rowbox{000}\rowbox{499}}}
&\hbox{\colbox{\rowbox{000}\rowbox{000}\rowbox{000}\rowbox{000}}}%
\hbox{\colbox{\rowbox{000}\rowbox{000}\rowbox{000}\rowbox{000}}}
&\hbox{\colbox{\rowbox{\textbf{77}}}}%
\hbox{\colbox{\rowbox{000}\rowbox{\textit{411}}\rowbox{422}\rowbox{433}}}%
\hbox{\colbox{\rowbox{444}\rowbox{\textit{455}}\rowbox{466}\rowbox{477}}}
&\hbox{\colbox{\rowbox{\textbf{27}}\rowbox{466}}}%
\hbox{\colbox{\rowbox{000}\rowbox{000}\rowbox{422}\rowbox{422}}}%
\hbox{\colbox{\rowbox{444}\rowbox{444}\rowbox{466}\rowbox{466}}}
\\
\hline
\hbox{\colbox{\hbox{5 }\vfill}}
&\hbox{\colbox{\rowbox{\textbf{77}}}}
&\hbox{\colbox{\rowbox{4EE}\rowbox{\textbf{77}}}}
&\hbox{\colbox{\rowbox{\textbf{57}}\rowbox{55C}\rowbox{5DD}}}
&\hbox{\colbox{\rowbox{444}\rowbox{5C5}\rowbox{444}\rowbox{5C5}}}
&\hbox{\colbox{\rowbox{4BB}}}%
\hbox{\colbox{\rowbox{11C}\rowbox{59D}\rowbox{\textbf{77}}\rowbox{\textbf{77}}}}
&\hbox{\colbox{\rowbox{4AA}\rowbox{4AA}}}%
\hbox{\colbox{\rowbox{000}\rowbox{\textit{585}}\rowbox{4AA}\rowbox{\textbf{77}}}}
&\hbox{\colbox{\rowbox{499}\rowbox{000}\rowbox{499}}}%
\hbox{\colbox{\rowbox{11C}\rowbox{\textbf{55}}\rowbox{11C}\rowbox{59D}}}
&\hbox{\colbox{\rowbox{000}\rowbox{000}\rowbox{000}\rowbox{000}}}%
\hbox{\colbox{\rowbox{000}\rowbox{\textit{585}}\rowbox{000}\rowbox{\textit{585}}}}
&\hbox{\colbox{\rowbox{\textit{\textbf{77}}}}}%
\hbox{\colbox{\rowbox{\textit{558}}\rowbox{\textit{559}}\rowbox{\textit{\textbf{77}}}\rowbox{\textbf{77}}}}%
\hbox{\colbox{\rowbox{55C}\rowbox{\textit{55D}}\rowbox{\textbf{77}}\rowbox{\textbf{77}}}}
&\hbox{\colbox{\rowbox{466}\rowbox{\textit{\textbf{67}}}}}%
\hbox{\colbox{\rowbox{000}\rowbox{\textit{141}}\rowbox{422}\rowbox{\textbf{77}}}}%
\hbox{\colbox{\rowbox{444}\rowbox{\textit{545}}\rowbox{466}\rowbox{\textbf{77}}}}
\\
\hline
\hbox{\colbox{\hbox{6 }\vfill}}
&\hbox{\colbox{\rowbox{\textbf{77}}}}
&\hbox{\colbox{\rowbox{\textbf{37}}\rowbox{6EE}}}
&\hbox{\colbox{\rowbox{4DD}\rowbox{6C6}\rowbox{\textbf{77}}}}
&\hbox{\colbox{\rowbox{444}\rowbox{444}\rowbox{6C6}\rowbox{6C6}}}
&\hbox{\colbox{\rowbox{4BB}}}%
\hbox{\colbox{\rowbox{22C}\rowbox{\textbf{77}}\rowbox{6AE}\rowbox{\textbf{77}}}}
&\hbox{\colbox{\rowbox{4AA}\rowbox{4AA}}}%
\hbox{\colbox{\rowbox{22C}\rowbox{22C}\rowbox{\textbf{33}}\rowbox{6AE}}}
&\hbox{\colbox{\rowbox{499}\rowbox{000}\rowbox{499}}}%
\hbox{\colbox{\rowbox{000}\rowbox{499}\rowbox{686}\rowbox{\textbf{77}}}}
&\hbox{\colbox{\rowbox{000}\rowbox{000}\rowbox{000}\rowbox{000}}}%
\hbox{\colbox{\rowbox{000}\rowbox{000}\rowbox{686}\rowbox{686}}}
&\hbox{\colbox{\rowbox{\textbf{77}}}}%
\hbox{\colbox{\rowbox{668}\rowbox{\textit{\textbf{77}}}\rowbox{66A}\rowbox{\textbf{77}}}}%
\hbox{\colbox{\rowbox{66C}\rowbox{\textit{\textbf{77}}}\rowbox{66E}\rowbox{\textbf{77}}}}
&\hbox{\colbox{\rowbox{\textbf{27}}\rowbox{666}}}%
\hbox{\colbox{\rowbox{\textbf{33}}\rowbox{668}\rowbox{\textbf{33}}\rowbox{66A}}}%
\hbox{\colbox{\rowbox{\textbf{37}}\rowbox{66C}\rowbox{\textbf{37}}\rowbox{66E}}}
\\
\hline
\hbox{\colbox{\hbox{7 }\vfill}}
&\hbox{\colbox{\rowbox{\textbf{77}}}}
&\hbox{\colbox{\rowbox{6EE}\rowbox{\textbf{77}}}}
&\hbox{\colbox{\rowbox{5DD}\rowbox{\textbf{77}}\rowbox{\textbf{77}}}}
&\hbox{\colbox{\rowbox{444}\rowbox{5C5}\rowbox{6C6}\rowbox{7C7}}}
&\hbox{\colbox{\rowbox{\textbf{77}}}}%
\hbox{\colbox{\rowbox{33C}\rowbox{\textbf{77}}\rowbox{\textbf{77}}\rowbox{\textbf{77}}}}
&\hbox{\colbox{\rowbox{4AA}\rowbox{\textbf{33}}}}%
\hbox{\colbox{\rowbox{22C}\rowbox{\textbf{77}}\rowbox{6AE}\rowbox{\textbf{77}}}}
&\hbox{\colbox{\rowbox{499}\rowbox{\textbf{11}}\rowbox{\textbf{55}}}}%
\hbox{\colbox{\rowbox{11C}\rowbox{59D}\rowbox{\textbf{77}}\rowbox{\textbf{77}}}}
&\hbox{\colbox{\rowbox{000}\rowbox{000}\rowbox{000}\rowbox{383}}}%
\hbox{\colbox{\rowbox{000}\rowbox{\textit{585}}\rowbox{686}\rowbox{787}}}
&\hbox{\colbox{\rowbox{\textbf{77}}}}%
\hbox{\colbox{\rowbox{\textbf{77}}\rowbox{\textit{\textbf{77}}}\rowbox{\textbf{77}}\rowbox{\textbf{77}}}}%
\hbox{\colbox{\rowbox{\textbf{77}}\rowbox{\textit{\textbf{77}}}\rowbox{\textbf{77}}\rowbox{\textbf{77}}}}
&\\
\hline
\hbox{\colbox{\hbox{8 }\vfill}}
&\hbox{\colbox{\rowbox{\textbf{77}}}}
&\hbox{\colbox{\rowbox{\textbf{33}}\rowbox{8EE}}}
&\hbox{\colbox{\rowbox{\textbf{55}}\rowbox{888}\rowbox{8DD}}}
&\hbox{\colbox{\rowbox{888}\rowbox{888}\rowbox{888}\rowbox{888}}}
&\hbox{\colbox{\rowbox{\textbf{11}}}}%
\hbox{\colbox{\rowbox{888}\rowbox{888}\rowbox{888}\rowbox{8BB}}}
&\hbox{\colbox{\rowbox{888}\rowbox{888}}}%
\hbox{\colbox{\rowbox{888}\rowbox{888}\rowbox{888}\rowbox{888}}}
&\hbox{\colbox{\rowbox{888}\rowbox{888}\rowbox{888}}}%
\hbox{\colbox{\rowbox{888}\rowbox{888}\rowbox{888}\rowbox{888}}}
&\hbox{\colbox{\rowbox{\textbf{01}}\rowbox{888}\rowbox{888}\rowbox{888}}}%
\hbox{\colbox{\rowbox{888}\rowbox{888}\rowbox{888}\rowbox{888}}}
&&\\
\hline
\hbox{\colbox{\hbox{9 }\vfill}}
&\hbox{\colbox{\rowbox{\textbf{77}}}}
&\hbox{\colbox{\rowbox{8EE}\rowbox{\textbf{77}}}}
&\hbox{\colbox{\rowbox{\textbf{55}}\rowbox{99C}\rowbox{9DD}}}
&\hbox{\colbox{\rowbox{888}\rowbox{9C9}\rowbox{888}\rowbox{9C9}}}
&\hbox{\colbox{\rowbox{8BB}}}%
\hbox{\colbox{\rowbox{99C}\rowbox{99C}\rowbox{99C}\rowbox{\textbf{55}}}}
&\hbox{\colbox{\rowbox{888}\rowbox{888}}}%
\hbox{\colbox{\rowbox{888}\rowbox{888}\rowbox{888}\rowbox{\textbf{11}}}}
&\hbox{\colbox{\rowbox{888}\rowbox{888}\rowbox{888}}}%
\hbox{\colbox{\rowbox{\textbf{55}}\rowbox{99C}\rowbox{99C}\rowbox{99C}}}
&&&\\
\hline
\hbox{\colbox{\hbox{A }\vfill}}
&\hbox{\colbox{\rowbox{\textbf{77}}}}
&\hbox{\colbox{\rowbox{\textbf{33}}\rowbox{AEE}}}
&\hbox{\colbox{\rowbox{8DD}\rowbox{ACA}\rowbox{\textbf{77}}}}
&\hbox{\colbox{\rowbox{888}\rowbox{888}\rowbox{ACA}\rowbox{ACA}}}
&\hbox{\colbox{\rowbox{8BB}}}%
\hbox{\colbox{\rowbox{AAC}\rowbox{AAC}\rowbox{AAC}\rowbox{\textbf{33}}}}
&\hbox{\colbox{\rowbox{888}\rowbox{888}}}%
\hbox{\colbox{\rowbox{\textbf{33}}\rowbox{AAC}\rowbox{AAC}\rowbox{AAC}}}
&&&&\\
\hline
\hbox{\colbox{\hbox{B }\vfill}}
&\hbox{\colbox{\rowbox{\textbf{77}}}}
&\hbox{\colbox{\rowbox{AEE}\rowbox{\textbf{77}}}}
&\hbox{\colbox{\rowbox{9DD}\rowbox{\textbf{77}}\rowbox{\textbf{77}}}}
&\hbox{\colbox{\rowbox{888}\rowbox{9C9}\rowbox{ACA}\rowbox{BCB}}}
&\hbox{\colbox{\rowbox{\textbf{11}}}}%
\hbox{\colbox{\rowbox{\textbf{77}}\rowbox{BBC}\rowbox{BBC}\rowbox{\textbf{77}}}}
&&&&&\\
\hline
\hbox{\colbox{\hbox{C }\vfill}}
&\hbox{\colbox{\rowbox{\textbf{77}}}}
&\hbox{\colbox{\rowbox{\textbf{37}}\rowbox{CEE}}}
&\hbox{\colbox{\rowbox{\textbf{57}}\rowbox{CCC}\rowbox{CDD}}}
&\hbox{\colbox{\rowbox{\textbf{07}}\rowbox{CCC}\rowbox{CCC}\rowbox{CCC}}}
&&&&&&\\
\hline
\hbox{\colbox[24pt]{\hbox{D }\vfill}}
&\hbox{\colbox[24pt]{\rowbox{\textbf{77}}}}
&\hbox{\colbox[24pt]{\rowbox{CEE}\rowbox{\textbf{77}}}}
&\hbox{\colbox[24pt]{\rowbox{\textbf{57}}\rowbox{DDC}\rowbox{DDD}}}
&&&&&&&\\
\hline
\hbox{\colbox[16pt]{\hbox{E }\vfill}}
&\hbox{\colbox[16pt]{\rowbox{\textbf{77}}}}
&\hbox{\colbox[16pt]{\rowbox{\textbf{37}}\rowbox{EEE}}}
&&&&&&&&\\
\hline
\hbox{\colbox[8pt]{\hbox{F }\vfill}}
&\hbox{\colbox[8pt]{\rowbox{\textbf{77}}}}
&&&&&&&&&\\
\hline
\end{tabular}

  \caption{General dimensions, part~1}
  \label{tab:n1}
\end{table}

\begin{table}[ptb]
  \centering
\newcommand{\tblfnt}{\fontsize{7}{8}\selectfont}%
\newcommand{\rowbox}[1]{\hbox to 18pt{\tblfnt \hfill\strut #1\hfill}}%
\newcommand{\colbox}[2][32pt]{\vbox to #1{\vfill #2}}%
\setlength{\tabcolsep}{1pt}%
\begin{tabular}{r|c|c|c|c|c|c|}
\setlength{\parindent}{0pt}
\hbox to 0pt{\hss$I\backslash$\hskip-3pt\raisebox{4pt}{$J$}}&5&4&3&2&1&0\\
\hline
\hbox{\colbox{\hbox{0 }\vfill}}
&\hbox{\colbox{\rowbox{\textit{\textbf{44}}}\rowbox{000}\rowbox{\textit{055}}}}%
\hbox{\colbox{\rowbox{000}\rowbox{000}\rowbox{000}\rowbox{000}}}%
\hbox{\colbox{\rowbox{000}\rowbox{\textit{055}}\rowbox{000}\rowbox{\textit{055}}}}
&\hbox{\colbox{\rowbox{000}\rowbox{000}\rowbox{000}\rowbox{000}}}%
\hbox{\colbox{\rowbox{000}\rowbox{000}\rowbox{000}\rowbox{000}}}%
\hbox{\colbox{\rowbox{000}\rowbox{000}\rowbox{000}\rowbox{000}}}
&\hbox{\colbox{\rowbox{000}}}%
\hbox{\colbox{\rowbox{000}\rowbox{000}\rowbox{000}\rowbox{000}}}%
\hbox{\colbox{\rowbox{000}\rowbox{000}\rowbox{000}\rowbox{000}}}%
\hbox{\colbox{\rowbox{000}\rowbox{000}\rowbox{000}\rowbox{000}}}
&\hbox{\colbox{\rowbox{000}\rowbox{000}}}%
\hbox{\colbox{\rowbox{000}\rowbox{000}\rowbox{000}\rowbox{000}}}%
\hbox{\colbox{\rowbox{000}\rowbox{000}\rowbox{000}\rowbox{000}}}%
\hbox{\colbox{\rowbox{000}\rowbox{000}\rowbox{000}\rowbox{000}}}
&\hbox{\colbox{\rowbox{000}\rowbox{000}\rowbox{000}}}%
\hbox{\colbox{\rowbox{000}\rowbox{000}\rowbox{000}\rowbox{000}}}%
\hbox{\colbox{\rowbox{000}\rowbox{000}\rowbox{000}\rowbox{000}}}%
\hbox{\colbox{\rowbox{000}\rowbox{000}\rowbox{000}\rowbox{000}}}
&\hbox{\colbox{\rowbox{\textbf{00}}\rowbox{000}\rowbox{000}\rowbox{000}}}%
\hbox{\colbox{\rowbox{000}\rowbox{000}\rowbox{000}\rowbox{000}}}%
\hbox{\colbox{\rowbox{000}\rowbox{000}\rowbox{000}\rowbox{000}}}%
\hbox{\colbox{\rowbox{000}\rowbox{000}\rowbox{000}\rowbox{000}}}
\\
\hline
\hbox{\colbox{\hbox{1 }\vfill}}
&\hbox{\colbox{\rowbox{\textit{\textbf{44}}}\rowbox{\textit{114}}\rowbox{\textit{155}}}}%
\hbox{\colbox{\rowbox{000}\rowbox{\textit{141}}\rowbox{000}\rowbox{\textit{141}}}}%
\hbox{\colbox{\rowbox{11C}\rowbox{\textit{\textbf{55}}}\rowbox{11C}\rowbox{\textit{15D}}}}
&\hbox{\colbox{\rowbox{000}\rowbox{\textit{141}}\rowbox{000}\rowbox{\textit{141}}}}%
\hbox{\colbox{\rowbox{000}\rowbox{\textit{141}}\rowbox{000}\rowbox{\textit{141}}}}%
\hbox{\colbox{\rowbox{000}\rowbox{\textit{141}}\rowbox{000}\rowbox{\textit{141}}}}
&\hbox{\colbox{\rowbox{000}}}%
\hbox{\colbox{\rowbox{\textit{114}}\rowbox{\textit{114}}\rowbox{\textit{114}}\rowbox{\textit{114}}}}%
\hbox{\colbox{\rowbox{000}\rowbox{000}\rowbox{000}\rowbox{000}}}%
\hbox{\colbox{\rowbox{11C}\rowbox{11C}\rowbox{11C}\rowbox{11C}}}
&\hbox{\colbox{\rowbox{000}\rowbox{000}}}%
\hbox{\colbox{\rowbox{000}\rowbox{000}\rowbox{000}\rowbox{000}}}%
\hbox{\colbox{\rowbox{000}\rowbox{000}\rowbox{000}\rowbox{000}}}%
\hbox{\colbox{\rowbox{000}\rowbox{000}\rowbox{000}\rowbox{000}}}
&\hbox{\colbox{\rowbox{000}\rowbox{000}\rowbox{000}}}%
\hbox{\colbox{\rowbox{\textit{\textbf{44}}}\rowbox{\textit{114}}\rowbox{\textit{114}}\rowbox{\textit{114}}}}%
\hbox{\colbox{\rowbox{000}\rowbox{000}\rowbox{000}\rowbox{000}}}%
\hbox{\colbox{\rowbox{\textbf{55}}\rowbox{11C}\rowbox{11C}\rowbox{11C}}}
&\\
\hline
\hbox{\colbox{\hbox{2 }\vfill}}
&\hbox{\colbox{\rowbox{\textit{055}}\rowbox{242}\rowbox{\textit{\textbf{66}}}}}%
\hbox{\colbox{\rowbox{000}\rowbox{000}\rowbox{242}\rowbox{242}}}%
\hbox{\colbox{\rowbox{000}\rowbox{\textit{055}}\rowbox{242}\rowbox{\textbf{77}}}}
&\hbox{\colbox{\rowbox{000}\rowbox{000}\rowbox{242}\rowbox{242}}}%
\hbox{\colbox{\rowbox{000}\rowbox{000}\rowbox{242}\rowbox{242}}}%
\hbox{\colbox{\rowbox{000}\rowbox{000}\rowbox{242}\rowbox{242}}}
&\hbox{\colbox{\rowbox{000}}}%
\hbox{\colbox{\rowbox{224}\rowbox{224}\rowbox{224}\rowbox{224}}}%
\hbox{\colbox{\rowbox{000}\rowbox{000}\rowbox{000}\rowbox{000}}}%
\hbox{\colbox{\rowbox{22C}\rowbox{22C}\rowbox{22C}\rowbox{22C}}}
&\hbox{\colbox{\rowbox{000}\rowbox{000}}}%
\hbox{\colbox{\rowbox{\textbf{22}}\rowbox{224}\rowbox{224}\rowbox{224}}}%
\hbox{\colbox{\rowbox{000}\rowbox{000}\rowbox{000}\rowbox{000}}}%
\hbox{\colbox{\rowbox{\textbf{33}}\rowbox{22C}\rowbox{22C}\rowbox{22C}}}
&&\\
\hline
\hbox{\colbox{\hbox{3 }\vfill}}
&\hbox{\colbox{\rowbox{\textit{155}}\rowbox{\textit{\textbf{66}}}\rowbox{\textit{\textbf{77}}}}}%
\hbox{\colbox{\rowbox{000}\rowbox{\textit{141}}\rowbox{242}\rowbox{343}}}%
\hbox{\colbox{\rowbox{11C}\rowbox{\textit{15D}}\rowbox{\textbf{77}}\rowbox{\textbf{77}}}}
&\hbox{\colbox{\rowbox{000}\rowbox{\textit{141}}\rowbox{242}\rowbox{343}}}%
\hbox{\colbox{\rowbox{000}\rowbox{\textit{141}}\rowbox{242}\rowbox{343}}}%
\hbox{\colbox{\rowbox{000}\rowbox{\textit{141}}\rowbox{242}\rowbox{343}}}
&\hbox{\colbox{\rowbox{000}}}%
\hbox{\colbox{\rowbox{\textbf{77}}\rowbox{334}\rowbox{334}\rowbox{334}}}%
\hbox{\colbox{\rowbox{\textbf{11}}\rowbox{338}\rowbox{338}\rowbox{338}}}%
\hbox{\colbox{\rowbox{\textbf{77}}\rowbox{33C}\rowbox{33C}\rowbox{33C}}}
&&&\\
\hline
\hbox{\colbox{\hbox{4 }\vfill}}
&\hbox{\colbox{\rowbox{\textit{\textbf{47}}}\rowbox{444}\rowbox{\textit{455}}}}%
\hbox{\colbox{\rowbox{000}\rowbox{\textit{411}}\rowbox{000}\rowbox{\textit{411}}}}%
\hbox{\colbox{\rowbox{444}\rowbox{\textit{455}}\rowbox{444}\rowbox{\textit{455}}}}
&\hbox{\colbox{\rowbox{\textbf{07}}\rowbox{444}\rowbox{444}\rowbox{444}}}%
\hbox{\colbox{\rowbox{000}\rowbox{000}\rowbox{000}\rowbox{000}}}%
\hbox{\colbox{\rowbox{444}\rowbox{444}\rowbox{444}\rowbox{444}}}
&&&&\\
\hline
\hbox{\colbox{\hbox{5 }\vfill}}
&\hbox{\colbox{\rowbox{\textit{\textbf{47}}}\rowbox{\textit{554}}\rowbox{\textit{555}}}}%
\hbox{\colbox{\rowbox{\textit{\textbf{55}}}\rowbox{\textit{\textbf{55}}}\rowbox{\textit{558}}\rowbox{\textit{559}}}}%
\hbox{\colbox{\rowbox{\textbf{57}}\rowbox{\textit{\textbf{57}}}\rowbox{55C}\rowbox{\textit{55D}}}}
&&&&&\\
\hline
\end{tabular}

  \caption{General dimensions, part~2}
  \label{tab:n2}
\end{table}

\begin{table}[ptb]
  \centering
  {\fontsize{10}{11}\selectfont
  \begin{tabular}[t]{|ccc|}
\hline
055&\textit{000}&\textbf{00}\\
057&\textit{000}&\textbf{00}\\
05D&\textit{000}&\textbf{00}\\
05F&\textit{000}&\textbf{00}\\
07D&\textit{000}&\textbf{00}\\
\hline
114&\textit{000}&\textbf{00}\\
115&\textit{000}&\textbf{00}\\
116&\textit{000}&\textbf{00}\\
117&\textit{000}&\textbf{00}\\
134&\textit{000}&\textbf{00}\\
135&\textit{000}&\textbf{00}\\
136&\textit{000}&\textbf{00}\\
137&\textit{000}&\textbf{00}\\
145&\textit{000}&\textbf{00}\\
147&\textit{000}&\textbf{00}\\
149&\textit{000}&\textbf{00}\\
14B&\textit{000}&\textbf{00}\\
14D&\textit{000}&\textbf{00}\\
14F&\textit{000}&\textbf{00}\\
155&\textit{000}&\textbf{00}\\
156&\textit{000}&\textbf{00}\\
\hline
\end{tabular}%
\hspace{10pt}\begin{tabular}[t]{|ccc|}\hline
157&\textit{000}&\textbf{00}\\
159&\textit{000}&\textbf{00}\\
15B&\textit{000}&\textbf{00}\\
15D&\textit{11C}&\textbf{11}\\
15F&\textit{11C}&\textbf{11}\\
167&\textit{066}&\textbf{22}\\
169&\textit{000}&\textbf{00}\\
16B&\textit{000}&\textbf{00}\\
16D&\textit{000}&\textbf{00}\\
177&\textit{077}&\textbf{33}\\
179&\textit{000}&\textbf{00}\\
17D&\textit{11C}&\textbf{11}\\
\hline
255&\textit{000}&\textbf{00}\\
257&\textit{242}&\textbf{22}\\
25D&\textit{000}&\textbf{00}\\
27D&\textit{22C}&\textbf{33}\\
\hline
345&\textit{000}&\textbf{00}\\
349&\textit{000}&\textbf{00}\\
34D&\textit{000}&\textbf{00}\\
355&\textit{000}&\textbf{00}\\
356&\textit{242}&\textbf{22}\\
\hline
\end{tabular}%
\hspace{10pt}\begin{tabular}[t]{|ccc|}\hline
357&\textit{343}&\textbf{33}\\
359&\textit{000}&\textbf{00}\\
35D&\textit{11C}&\textbf{11}\\
369&\textit{000}&\textbf{00}\\
36D&\textit{22C}&\textbf{33}\\
379&\textit{338}&\textbf{11}\\
37D&\textit{33C}&\textbf{33}\\
\hline
455&\textit{444}&\textbf{03}\\
457&\textit{444}&\textbf{03}\\
459&\textit{000}&\textbf{00}\\
45B&\textit{000}&\textbf{00}\\
45D&\textit{444}&\textbf{03}\\
45F&\textit{444}&\textbf{03}\\
479&\textit{000}&\textbf{00}\\
47D&\textit{444}&\textbf{03}\\
\hline
555&\textit{444}&\textbf{03}\\
556&\textit{444}&\textbf{03}\\
557&\textit{444}&\textbf{03}\\
558&558&\textbf{\textit{10}}\\
559&\textit{558}&\textbf{\textit{10}}\\
55A&558&\textbf{\textit{10}}\\
\hline
\end{tabular}%
\hspace{10pt}\begin{tabular}[t]{|ccc|}\hline
55B&\textit{558}&\textbf{\textit{10}}\\
55D&\textit{55C}&\textbf{13}\\
55F&\textit{55C}&\textbf{13}\\
567&\textit{466}&\textbf{23}\\
569&\textit{000}&\textbf{00}\\
56D&\textit{444}&\textbf{03}\\
577&\textit{477}&\textbf{33}\\
578&558&\textbf{\textit{10}}\\
579&\textit{558}&\textbf{\textit{10}}\\
57A&57A&\textbf{\textit{32}}\\
57D&\textit{55C}&\textbf{13}\\
58D&585&\textbf{\textit{10}}\\
58F&585&\textbf{\textit{10}}\\
5AD&585&\textbf{\textit{10}}\\
\hline
679&\textit{668}&\textbf{33}\\
67D&\textit{66C}&\textbf{33}\\
\hline
779&\textit{778}&\textbf{33}\\
77D&\textit{77C}&\textbf{33}\\
78D&585&\textbf{\textit{10}}\\
\hline
\end{tabular}}

  \caption{Differences in dimension $8$}
  \label{tab:d8}
\end{table}

\begin{table}[tbp]
  \centering
\newcommand{\rowbox}[1]{\hbox to 32pt{\hfill\strut #1\hfill}}%
\newcommand{\colbox}[2][26pt]{\vbox to #1{\vfill #2}}%
\setlength{\tabcolsep}{1pt}%
\begin{tabular}{r|c|c|c|c|}
\setlength{\parindent}{0pt}
\hbox to 0pt{\hss$I\backslash$\hskip-3pt\raisebox{4pt}{$J$}}&A&8&2&0\\
\hline
\hbox{\colbox{\hbox{\strut 0 }\vfill}}
&\hbox{\colbox{\rowbox{\textbf{0AA}}}}
&\hbox{\colbox{\rowbox{000}\rowbox{000}}}
&\hbox{\colbox{\rowbox{000}}}%
\hbox{\colbox{\rowbox{000}\rowbox{000}}}
&\hbox{\colbox{\rowbox{\textbf{000}}\rowbox{000}}}%
\hbox{\colbox{\rowbox{000}\rowbox{000}}}
\\
\hline
\hbox{\colbox{\hbox{\strut 2 }\vfill}}
&\hbox{\colbox{\rowbox{\textbf{2AA}}}}
&\hbox{\colbox{\rowbox{000}\rowbox{282}}}
&\hbox{\colbox{\rowbox{000}}}%
\hbox{\colbox{\rowbox{\textbf{228}}\rowbox{228}}}
&\\
\hline
\hbox{\colbox{\hbox{\strut 8 }\vfill}}
&\hbox{\colbox{\rowbox{\textbf{8AA}}}}
&\hbox{\colbox{\rowbox{\textbf{888}}\rowbox{888}}}
&&\\
\hline
\hbox{\colbox[20pt]{\hbox{\strut A }\vfill}}
&\hbox{\colbox[20pt]{\rowbox{\textbf{AAA}}}}
&&&\\
\hline
\end{tabular}

  \caption{Dimension $4$}
  \label{tab:d4}
\end{table}

\addcontentsline{toc}{section}{References}

{\small
\noindent
(Mart\'\i n Cabrera) Department of Fundamental Mathematics, University of
La Laguna, 38200 La Laguna, Tenerife, Spain

{\it E-mail address}: {\tt fmartin@ull.es}
\smallskip

\noindent
(Swann) Department of Mathematics and Computer Science, University of
Southern Denmark, Campusvej 55, DK-5230 Odense M, Denmark

{\it E-mail address}: {\tt swann@imada.sdu.dk}
}
\end{document}